\documentclass[12pt,a4paper,final]{iopart}

\pdfoutput=1
\usepackage{iopams}  
\usepackage{algorithm}
 \usepackage{algorithmicx}
 \usepackage{bbm}

\newtheorem{Remark}{Remark}[section]

\newcommand{\bbR}{\mathbb{R}}
\newcommand{\bbP}{\mathbb{P}}

\usepackage{graphicx}
\usepackage{color}

\expandafter\let\csname equation*\endcsname\relax

\expandafter\let\csname endequation*\endcsname\relax

\usepackage{amsmath}

\begin{document}

\title[Parameterizations for Ensemble Kalman Inversion]{Parameterizations for Ensemble Kalman Inversion}

\author{Neil K. Chada$^{*}$, Marco A. Iglesias$^{**}$, Lassi Roininen$^{\dagger}$ \\ and Andrew M. Stuart$^{\ddagger}$}

\address{$^*$ Mathematics Institute, University of Warwick, Coventry, CV4 7AL, UK. \\
$^{**}$ School of Mathematical Sciences, University of Nottingham, \\ Nottingham Park, NG7 2QL, UK. \\
$^{\dagger}$ Department of Mathematics, Imperial College London, London, SW7 2RD, UK. \\
$^{\ddagger}$ Computing and Mathematics Sciences, California Institute of Technology, \\ Pasadena, CA 91125, USA.}
\ead{\\ n.chada@warwick.ac.uk \\
marco.iglesias@nottingham.ac.uk \\
l.roininen@imperial.ac.uk \\
astuart@caltech.edu}

\begin{abstract} The use of ensemble methods to solve inverse problems is attractive 
because it is a derivative-free methodology which is also
well-adapted to parallelization. 
In its basic iterative form the method produces an ensemble of solutions which lie
in the linear span of the initial ensemble. Choice of the parameterization of the unknown
field is thus a key component of the success of the method. We demonstrate how both
geometric ideas and hierarchical ideas can be used to design effective parameterizations
for a number of applied inverse problems arising in electrical impedance tomography, groundwater flow and source inversion. In particular we show how geometric ideas, including the level set method,
can be used to reconstruct piecewise continuous fields, and we show how hierarchical
methods can be used to learn key parameters in continuous fields, such as length-scales,
resulting in improved reconstructions. Geometric and hierarchical ideas are combined
in the level set method to find piecewise constant reconstructions with interfaces
of unknown topology.
\end{abstract}

\vspace{2pc}
\submitto{Inverse Problems}

\newpage

\section{Introduction}
\subsection{Content}
Consider finding $u$ from $y$ where
\begin{equation}
\label{inverse}
y=\mathcal{G}(u)+\eta,
\end{equation}
$\mathcal{G}$ is a forward map taking 
the unknown parameter $u$ into the data space,  
and $\eta$ represents noise. 
Ensemble Kalman inversion is an attractive technique that has shown considerable success in the solution of such problems. 
Whilst it is derived from the application of Kalman-like thinking, with means and covariances computed from an empirical ensemble, it essentially acts as a black-box derivative-free optimizer which requires only evaluation of the forward map $\mathcal{G}(\cdot);$ in practice it can often return good solutions to inverse problems with relatively few forward map
evaluations. However the choice of parameterization of the unknown is key to the
success of the method. In this paper we will demonstrate how carefully thought out
parameterizations can have substantial impact in the quality of the reconstruction. 

Although our viewpoint in this paper is to consider Ensemble Kalman 
inversion as an optimization method, 
and evaluate it from this perspective, there is considerable insight to be gained
from the perspective of Bayesian inversion; this is despite the fact
that the the algorithm does not, in general,
recover or sample the true Bayesian posterior distribution of the 
inverse problem. Algorithms that can, with controllable error, approximately sample from the true posterior distribution are commonly referred to as fully Bayesian, with examples including Markov Chain Monte Carlo and sequential Monte Carlo.
Ensemble Kalman inversion is not fully Bayesian but the link to 
Bayesian inversion remains important as we now explain.  
There is considerable literature available about methods to improve fully Bayesian approaches to the inverse problem through for example, geometric and hierarchical 
parameterizations of the unknown.
The purpose of this paper is to demonstrate how these ideas from Bayesian inversion may be used with some success to improve the capability of ensemble Kalman methods considered as optimizers. In view of the relatively low computational cost of the ensemble methods, in comparison with fully Bayesian inversion, this cross-fertilization of ideas has the potential to be quite fruitful.

\subsection{Literature Review}

The Kalman filter (KF) \cite{Kalman} was 
developed to sequentially update the probability distribution on
states of partially observed linear Gaussian systems, 
and subsequently generalized to nonlinear problems in the form of the extended
Kalman filter. However for high-dimensional systems  the size of
covariances makes use of these methods prohibitive. In 1994 Evensen \cite{OG,EnKF} proposed a 
Monte-Carlo based nonlinear Kalman filter which tackled this issue by using an ensemble 
of particles to represent the covariances and mean, resulting in what is now
known as the ensemble Kalman filter (EnKF).
A major success story for the EnKF has been in weather prediction models 
\cite{Anderson,Weather}, but it has also been deployed in numerous applications domains,
including  the reservoir engineering community \cite{Aan} and in oceanography 
\cite{VL}.
Variants on the idea include the randomized maximum likelihood (RML) method
\cite{Oliver}, 
and algorithms such as the ensemble square-root Kalman filter \cite{ESRKF}.

In this paper we are primarily interested in use of ensemble Kalman methods to
study inverse problems for parameter estimation, an approach pioneered for 
oil industry applications where the inverse problem is known as
history matching \cite{Li, Oliver}; the 
paper \cite{ernst2015analysis} contains an insightful analysis of the 
methodology in the large ensemble limit. 
In this application domain such inversion 
methods are sometimes referred to as ensemble Kalman smoothers, although 
the nomenclature is not uniform. We will simply refer to 
{\em ensemble Kalman inversion} (EKI).
The methodology is formulated quite generally in \cite{ILS2}, independently
of oil industry applications; and in \cite{ILS} it is shown that the method
performs well as an optimizer but does not capture true posterior
uncertainty, in the context of  oil industry applications.

The ideas introduced in this paper 
concerning parameterization are independent of 
the particular implementation of the EKI method used;
in all our numerical experiments we will use the form of iterative 
regularization proposed by Iglesias \cite{reg}. 
The general philosophy behind the method is that, as the algorithm is iterated, 
the solution to the inverse problem should approach the truth in the small noise
limit, and hence that any regularization introduced should diminish in
influence in this limit. However, because the convergence theory for ensemble Kalman
inversion is in its infancy, the choice of iterative regularization is made
by analogy with  classical  iterative methods that have been used for inverse 
problems \cite{ IM, GN, Hanke}, with the ensemble method using empirical covariances
in place of derivatives or adjoints of the forward solver. 
The resulting iterative method is an ensemble version of the Levenberg-Marquardt 
algorithm with the inclusion of regularization as in \cite{Hanke}.

In the form of EKI that we use the linear span of
the initial ensemble is preserved by the iteration \cite{ILS2,Li}.
The initial ensemble thus encodes prior information about the solution of the problem.
This means that choice of the parameterization of the method, as well as the
choice initial subspace, is key to its performance.
Based on experience with (Bayesian) statistical modelling we will introduce
geometric and hierarchical priors that address the issue of making good
parameterizations, and we will draw from those priors to create the initial
ensemble.
Hierarchical models have been extensively studied in the fields of computational statistics and machine learning \cite{Liang,centre1,Teh}. Their use
in the context of Monte Carlo Markov chain (MCMC) methods for
Bayesian inverse problems is overviewed in \cite{centre2}; see 
\cite{MD2, raul} for application oriented work. One important outcome
of research in this area is that
learning parameters such as length-scale, amplitude and regularity
within Gaussian random field priors (such as {Whittle-Mat\'{e}rn}) 
can be of significant value \cite{Sergios,Liu,lassi}.
In a series of recent papers this hierarchical modelling was extended to allow
for length-scale which is itself spatially varying \cite{Markku, lassi}. 
The development of hierarchical methods within EKI, rather than fully
Bayesian MCMC, has been
limited to date, with the primary contribution being the work
\cite{Hier} where the methodology was based on building large ensembles 
from multiple Gaussians assigned different weights.  
However this work requires that the correct hierarchical 
parameter is in the ensemble if it is to be successful. 
We also note that there is some work in hierarchical
EnKF within the context of state estimation; see
\cite{SPE-2} and the {references} therein.

In addition to hierarchical approaches we will also study geometric
parameterizations. These can be of use when the geometric object has known
form, such as faults, channels and inclusions, or when it is
of unknown topology the level set method may be used \cite{Santosa1996}.
We will build on recent Bayesian implementation of these ideas; see
\cite{geo,Lu} and references therein.

\subsection{Contribution of This Work}
Our main contribution is to establish the importance of
novel parameterizations which
have the potential to substantially improve the performance of EKI. 
Although our perspective on EKI is one of optimization, the methods we
introduce are all based on taking established and emerging methods 
from Bayesian statistics and developing them in the context of the 
ensemble methods. The connection to Bayesian statistics is exploited 
to provide insights about how to make these methods efficient.  
The resulting methods are illustrated by means of examples arising 
in both electrical impedance tomography, groundwater flow and source
inversion. The contributions are: \smallskip
 \begin{itemize}
 \item We develop hierarchical approaches for EKI, based on solving for the unknown function and unknown scalars which parameterize the prior.
 \item We generalize these hierarchical approaches to EKI
to include unknown fields which parameterize the prior, rather than scalars. 
 \item We demonstrate the key role of choosing non-centered variables when implementing hierarchical methods. 
 \item  We show the potential for geometric hierarchical priors, including the level set parameterization, for piecewise continuous reconstructions.
 \end{itemize}

\subsection{Organization}
The layout of the paper is as follows, in section \ref{section-2} we discuss different
approaches to parameterizing inverse problems. We begin by conveying the main ideas in 
subsection \ref{subsection-2.1} in abstract. In subsection \ref{subsection-2.2} we describe 
these ideas more concretely. 
In section \ref{section-2add} we describe the hierarchical version of iterative 
EKI as used in this paper, and section \ref{section_model} describes
the model problems that we use to 
illustrate the power of the proposed parameterizations. Numerical results
are presented in section \ref{section-4}, {whilst in section 
\ref{section-5} we make some concluding remarks.}
\subsection{Notation}
Throughout the paper we make use of common notation for Hilbert space norms and inner products, $\| \cdot \|, \langle \cdot \rangle$. We will assume that $\mathcal{X}$ and $\mathcal{Y}$ are two separable Hilbert spaces which are linked through the forward operator $\mathcal{G}: \mathcal{X} \rightarrow \mathcal{Y}$. This nonlinear operator can be thought of as mapping from the space of {unknown parameters} $\mathcal{X}$ to the observation space $\mathcal{Y}$. Our additive noise for the inverse problems will be denoted by  $\eta \sim N(0,\Gamma)$ where $\Gamma: \mathcal{Y} \rightarrow \mathcal{Y}$ is a self-adjoint positive operator. For any such operator we define $ \langle \cdot, \cdot \rangle_{\Gamma} = \langle \Gamma^{-1/2}\cdot, \Gamma^{-1/2}\cdot \rangle$ and $\| \cdot \|_{\Gamma} = \| \Gamma^{-1/2} \cdot \|$, and for finite dimensions $| \cdot |_{\Gamma} = | \Gamma^{-1/2} \cdot|$ with $|\cdot|$the Euclidean norm. {If the Gaussian measure associated to $\eta$ is supported on ${\mathcal Y}$ then we will require $\Gamma$ to be trace-class; however we will also consider white noise whose support is on a larger space than  ${\mathcal Y}$ and for which the trace-class condition fails in the infinite dimensional setting.}

\section{Inverse Problem} \label{section-2}
\subsection{Main Idea} \label{subsection-2.1}
\subsubsection{Non-Hierarchical Inverse Problem} \label{subsubsection-2.1.1}

We are interested in the recovery of $ u \in \mathcal{X}$ from measurements of  $y \in \mathcal{Y}$ given by equation \eqref{inverse}
in which, recall, $\eta$ is additive Gaussian noise. In the Bayesian approach to inverse problems we treat each quantity within (\ref{inverse}) as a random variable. Via an application of Bayes' Theorem \footnote{{We write all instances of Bayes' Theorem in finite dimensions for simplicity; extension to Bayes' Theorem for functions
is straightforward but not central to this paper and so we avoid the extra notation that would be needed for this.}} 
{\cite{BIP}} we can characterize the conditional distribution of $u|y$ via
\begin{equation}
\label{b1}
\mathbb{P}(u|y) \propto \mathbb{P}({y|u}) \times \mathbb{P}({u}),
\end{equation}
where $\mathbb{P}({u})$ is the prior distribution,
$\mathbb{P}(u|y)$ is the posterior distribution and $\mathbb{P}({y|u})$
is the likelihood. 
Although we view EKI as an optimizer in this paper,
the Bayesian formulation of \eqref{inverse} is important because we derive
the initial ensemble from the prior distribution $\mathbb{P}({u}).$ 
From an optimization viewpoint, our goal is to make the following least
squares objective function small:

\begin{equation}
\label{eq:fi}
\Phi(u;y)=\frac12 |y-\mathcal{G}(u)|^2_{\Gamma}.
\end{equation}
This (upto an irrelevant additive constant)
is the negative log likelihood since 
it is assumed that $\eta$ is a mean-zero
Gaussian with covariance $\Gamma.$

\subsubsection{Centered Hierarchical Inverse Problem}
\label{subsubsection-2.1.2}

In many applications it can be advantageous to add additional unknowns $\theta$ 
to the inversion process. In particular these may enter through the prior, 
as in hierarchical methods, \cite{centre2}, and we will refer to such
parameters as hyperparameters. 
The inverse problem is then the recovery of $(u, \theta)$ from measurements of $y$ given again by \eqref{inverse}. 
The additional parameterization of the prior results in Bayes' Theorem in the form
\begin{equation}
\label{b2}
\mathbb{P}(u,\theta|y) \propto \mathbb{P}(y|u) \times \mathbb{P}({u, \theta}).
\end{equation}
Prior samples, used to initialize the ensemble smoother,  will then be of the pair
$(u, \theta).$
What we will term {\em centered} hierarchical methods (a terminology we discuss
in Section \ref{section-2add}) typically involve factorization
of the prior in the form $\mathbb{P}({u, \theta})=\mathbb{P}(u|\theta)\mathbb{P}(\theta).$ From an optimization point of view our goal is
again to make the objective function \eqref{eq:fi} small, but now
using hierarchical parameterization to construct the initial ensemble. \footnote{We note that hyperparameters $\theta$
may also enter the likelihood as well as the prior if the state variable
is re-scaled in a hyperparameter-dependent fashion, as happens in the 
version of the Bayesian level set method advocated in \cite{MD2}.}

\subsubsection{Non-Centered Hierarchical Inverse Problem}
\label{subsubsection-2.1.3}

Another variant of the inverse problem that is particularly relevant for hierarchical methods, in which $\theta$ enters only the prior, is {\em non-centered} reparameterization
(a further terminology discussed in Section \ref{section-2add}). 
We introduce the transformation $T:(\xi,\theta) \rightarrow u$ 
and note that (\ref{inverse}) then becomes
\begin{equation}
\label{hier_inv2}
y = \mathcal{G}(T(\xi, \theta)) + \eta,
\end{equation}
and Bayes' Theorem then reads 
\begin{equation}
\label{b3}
\mathbb{P}(\xi, \theta|y) \propto \mathbb{P}({y|\xi, \theta}) \times \mathbb{P}({\xi, \theta}).
\end{equation}
Prior samples, again used to initialize the ensemble smoother,  will then be of the pair
$(\xi, \theta).$
Typically the change of variables from $u$ to $\xi$ is introduced so that
$\xi$ and $\theta$ are independent under the prior: 
$\mathbb{P}({\xi, \theta})=\mathbb{P}(\xi)\mathbb{P}(\theta).$ {As a result the inverse problem (\ref{hier_inv2}) is different to that 
appearing in (\ref{inverse}), in terms of both the prior and the likelihood. This non-centered approach is equivalent to the ancillary augmentation technique discussed in \cite{centre1} which discusses the decoupling of $u$ and $\theta$ through the variable $\xi$.}
From an optimization viewpoint, our goal is to make the following least
squares objective function small:
\begin{equation}
\label{eq:fi3}
\Phi(\xi,\theta;y)=\frac12 |y-\mathcal{G}(T(\xi,\theta))|^2_{\Gamma}.
\end{equation}

{
In this paper we consider the application of EKI for
solution of the inverse problems (\ref{inverse}), non-hierarchically
and centered-hierarchically, and (\ref{hier_inv2}) non-centered hierarchically.
Although this method has a statistical derivation, the work in \cite{ILS2, LS} demonstrates that the method may be thought of as a derivative-free optimizer that approximates the least squares problem (\ref{eq:fi}) and (\ref{eq:fi3}) rather than sampling from the relevant Bayesian posterior
distribution. We will show that the use of iterative EKI, 
as proposed by Iglesias in \cite{reg}, can effectively solve a wide range of challenging 
inversion problems, 
if judiciously parameterized.
}
\\

\subsection{Details of Parameterizations}  \label{subsection-2.2}
In this section we describe in detail several classes of parameterizations that
we use in this paper. The first and second are geometric parameterizations,
ideal for piecewise continuous reconstructions with unknown interfaces. The third
and fourth are hierarchical methods which introduce an unknown length-scale,
and regularity parameter, into the inversion.
For the two geometric problems we initially formulate 
in terms of trying to find a function $w: D \mapsto \bbR$, $D$
a subset of $\bbR^d$, and then reparameterize $w.$
For the hierarchical problems  we initially formulate 
in terms of trying to find a function $u: D \mapsto \bbR$, and
then append parameters $\theta$ and also rewrite in terms of $(\xi,\theta)
\mapsto u.$

\subsubsection{Geometric Approach -- Finite Dimensional Parameterization}
\label{sssec:g}
In many problems of interest the unknown function $w$ has discontinuities, determination of which forms part of the
solution of the inverse problem. To tackle such problems it may be useful to write
$w$ in the form
\begin{equation}
w(x) = \sum^{n}_{i=1}u_{i}(x)\chi_{D_i}(x).
\end{equation}
Here the union of the disjoint sets $D_i$ is the whole domain $D$. If we assume that the 
configuration of the $D_i$ is determined by a finite
set of scalars $\theta$ and let $u$ denote the union of the functions $u_i$ and
the parameters $\theta$ 
then we may rewrite the inverse problem in the form \eqref{inverse}. 
{The case where the number of subdomains $n$ is unknown would be an interesting and useful extension of this work; but we do not consider it here.}

\subsubsection{Geometric Approach -- Infinite Dimensional Parameterization}
\label{sssec:l}
If the interface boundary is not readily described by a finite number of parameters we may use the level set idea. For example if the field $w$ takes two 
known values $w^{\pm}$ 
with unknown interfaces between them we may write
\begin{equation}
w(x)=w^{+}{\mathbb I}_{u>0}(x)+w^{-}{\mathbb I}_{u<0}(x),
\end{equation}
and formulate the inverse problem in the form \eqref{inverse} for $u$.
This idea may be generalized to functions which take an arbitrary number
of constant values, through the introduction of level sets other than $u=0$,
or through vector level sets functions $u.$
\\
\subsubsection{Scalar-valued Hierarchical Parameterizations}
\label{sssec:h}
To illustrate ideas we will concentrate on Gaussian priors of Whittle-Mat\'{e}rn type. These are characterized by a covariance function of the form
\begin{equation}
\label{fun}
c(x,y) = \sigma^2 \frac{2^{1-\alpha + d/2}}{\Gamma(\alpha - d/2)}\bigg(\frac{|x-y|}{\ell}\bigg)^{\alpha - d/2} {K_{\alpha - d/2}\bigg(\frac{|x-y|}{\ell}\bigg)},  \ \ \ x,y \in \mathbb{R}^d,
\end{equation}
where $K_{\cdot}$ is a modified Bessel function of the second kind, $\sigma^2>0$ is the variance and $\Gamma(\cdot)$ is a Gamma function. {We will always ensure that $\ell>0$ and $\alpha>d/2$ so that draws from the Gaussian are well-defined and continuous.} On the unbounded domain $\bbR^d$, samples
from this process may be generated by solving the stochastic PDE 
\begin{equation}
\label{SPDE0}
(I-\ell^2 \triangle)^{\frac{\alpha}{2}}u = \ell^{d/2}\sqrt{\beta}\xi,
\end{equation}
where $\xi  \in H^{-s}(D)$, $s>\frac{d}{2}$, is a Gaussian white noise,
i.e. $ \xi \sim N (0,I)$, and 
$$\beta=\sigma^2\frac{2^d\pi^{d/2}\Gamma(\alpha)}{\Gamma(\alpha-\frac{d}{2})}.$$

In this paper we will work with the scalar hierarchical parameters $\alpha$
and $\tau=\ell^{-1}$. Putting $\ell=\tau^{-1}$ into \eqref{SPDE0} 
gives the stochastic PDE
\begin{equation}
\label{SPDE}
\mathcal{C}_{\alpha,\tau}^{-\frac{1}{2}}u = \xi,
\end{equation}
where the covariance operator $\mathcal{C}_{\alpha,\tau}$ has the form 
$\mathcal{C}_{\alpha,\tau} = \tau^{2\alpha-d} \beta(\tau^2 I-  \Delta )^{-\alpha}.$
Throughout this paper we choose $\beta=\tau^{d-2\alpha}$ so that
\begin{equation}
\label{cov_prior}
\mathcal{C}_{\alpha,\tau} =(\tau^2 I-  \Delta )^{-\alpha},
\end{equation}
and we equip the operator $\Delta$ with Dirichlet boundary conditions
on $D$. {These} choices simplify the exposition, but are not an integral
part of the methodology; different choices could be made.

\begin{figure}[h!]
\centering
 \includegraphics[scale=0.43]{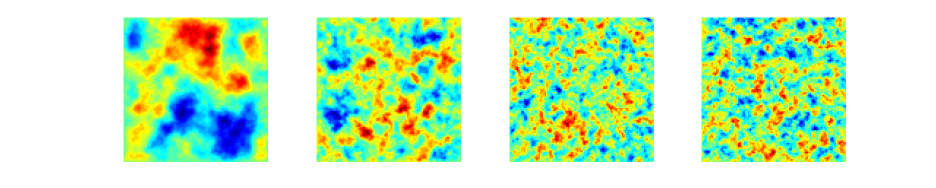}
\caption{Modified inverse length scale for $\tau$ = 10, 25, 50 and 100. Here $\alpha=1.6$.}
 \label{matern}
\end{figure}

\begin{figure}[h!]
\centering
 \includegraphics[scale=0.55]{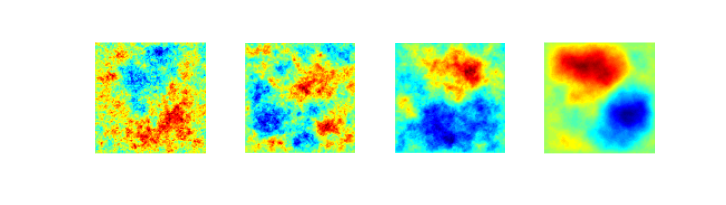}
\caption{Modified regularity for $\alpha$ = 1.1, 1.3, 1.5 and 1.9. Here $\tau=15$.}
 \label{alpha}
\end{figure}

{We have thus formulated the inverse
problem in the form of the centered hierarchical inverse problem
\eqref{inverse} for $(u,\theta)$ with 
$\theta=(\alpha,\tau)$. The parameter $\theta$ enters only through the prior as in this
particular case it does not appear in the likelihood.} 
In this paper we will place uniform priors on $\alpha$ and $\tau$, {which will be specified in subsection \ref{subsection-4.1}.}
We may also work with the variables $(\xi,\theta)$ noting that \eqref{SPDE} 
defines a map $T: (\xi,\theta) \mapsto u$ and we have formulated the inverse
problem in the form \eqref{hier_inv2}, 
{the non-centered hierachical form.} 
In Figures \ref{matern} and \ref{alpha} we display random samples from (\ref{cov_prior}) with imposed Dirichlet boundary conditions and varying values of the inverse length scale $\tau$ and the regularity $\alpha$. {These samples are constructed in the domain $D=[0,1]^2$.}

\subsubsection{Function-Valued Hierarchical Parameterization}\label{subsubsection-2.2.3}
\begin{figure}[h!]
\centering
 \includegraphics[scale=0.45]{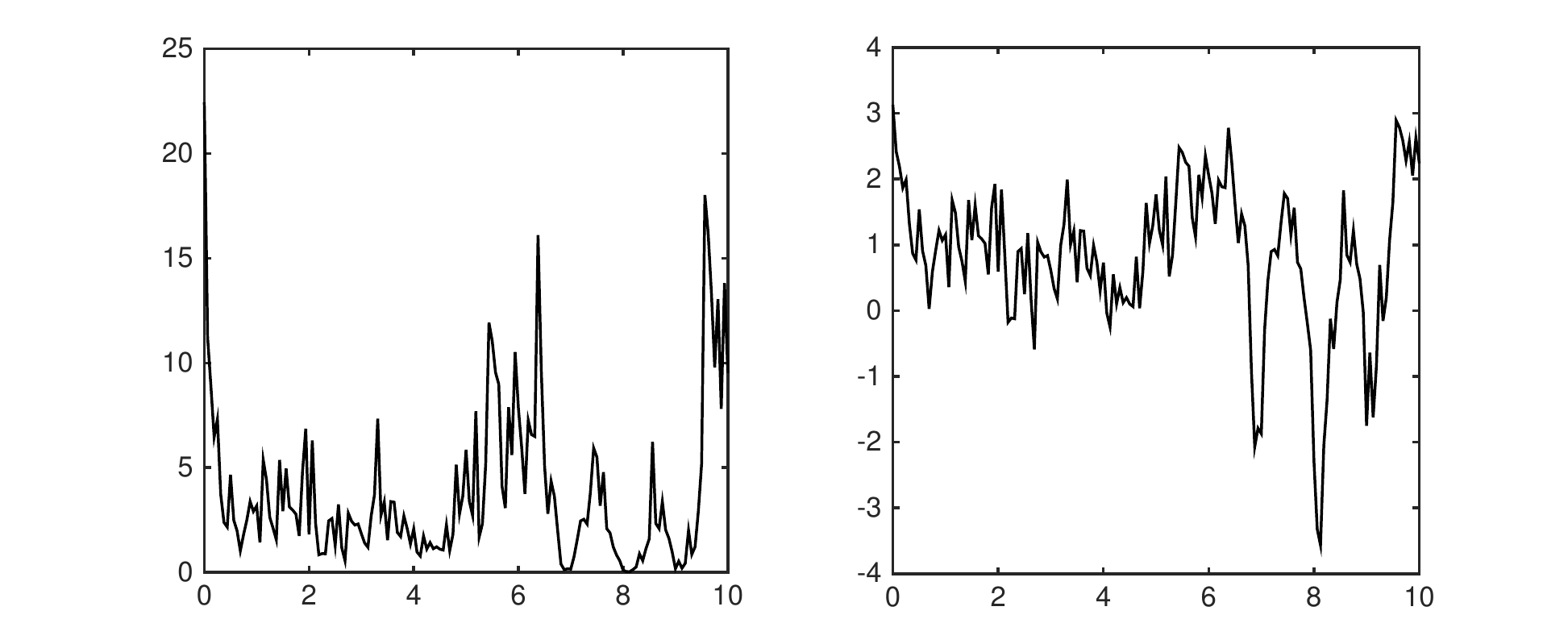}
\caption{Gaussian random field. Left: Length-scale realization $\ell(x)$. Right: Realization of $v(x)$. }
 \label{fig:gauss_real_walk}
\end{figure}

\begin{figure}[h!]
\centering
 \includegraphics[scale=0.45]{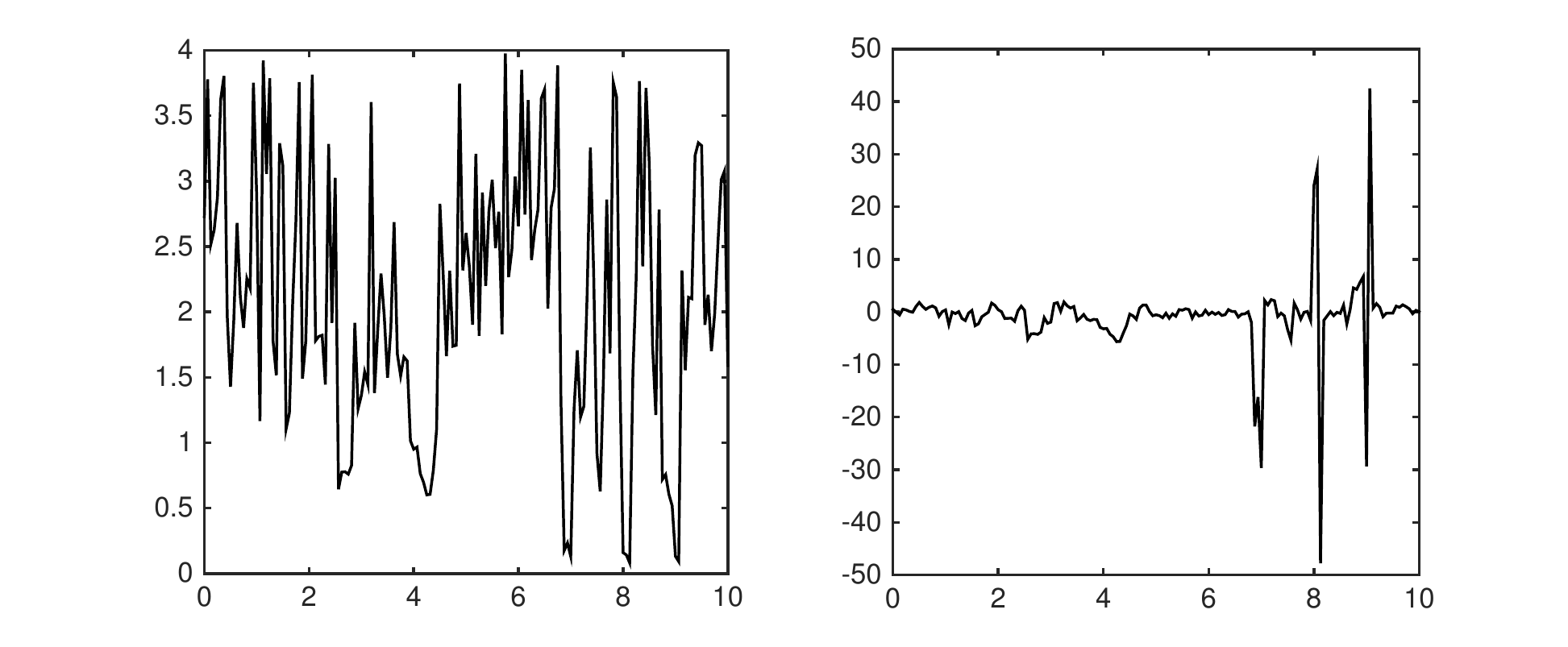}
\caption{Cauchy random field. Left: Length-scale realization $\ell(x)$. Right: Realization of $v(x)$. }
 \label{fig:gauss_real_length}
\end{figure}

In order to represent non-stationary features it is of interest to allow hierarchical parameters to themselves vary in space. To this end we will also seek to generalize \eqref{SPDE0} and work with the form
\begin{equation}
\label{SPDE2}
\left(I -\ell(x;v)^2 \Delta\right)^{\frac{\alpha}{2}}u = \ell(x;v)^{\frac{d}{2}} \xi.
\end{equation}
(We have set $\beta=1$ for simplicity).
In order to ensure that the length-scale is positive we will write it in the form
\begin{equation}
\label{eq:gp}
\ell(x;v) = g(v(x)),
\end{equation}
for some {positive} monotonic increasing function $g(\cdot)$.
We thus have a formulation as a {centered hierarchical
inverse problem in the form \eqref{inverse} 
noting that hyperparameter  $\theta=v$ is here a function and 
enters only through the prior.} 
We may also formulate inversion in terms of
the variables $(\xi,v)$ giving the inverse problem in the form \eqref{hier_inv2} 
with $\theta=v$. We will consider two forms of prior on $v$. The first is 
based on a Gaussian random field with Whittle-Mat\'{e}rn covariance function \eqref{fun}
and we then choose $g(v)=\exp(v).$ 
The second, which will apply only in one dimension, is to consider 
a one-dimensional Cauchy process, as in \cite{Markku}.
{In particular  we will construct $\ell(x;v)$ by employing
a one-dimensional Cauchy process $v(x)$, which is an
$\alpha$-stable L\'{e}vy motion with $\alpha =1$ with Cauchy increments
on the interval $\delta$ given by the density function $f$, and positivity-inducing
function $g$, where}
\begin{equation}
\label{eq:cp}
f(x) = \frac{\delta}{\pi(\delta^2 + x^2)}, \quad g(s) = \frac{a}{b+c|s|}+d,
\end{equation}
{such that} $a,b,c,d>0$ are constants. \footnote{{We note here that $\alpha$ has a different meaning from the parameter
$\alpha$ used in the covariance function of a Gaussian prior in, for example,
\eqref{cov_prior}. We abuse notation in this way because the parameter
$\alpha$ is widely used in the literature in both contexts.}}
Samples from these two priors on $v$, and hence $\ell$, are 
shown in Figures \ref{fig:gauss_real_walk} and \ref{fig:gauss_real_length}. 

\section{Iterative Ensemble Kalman Inversion} \label{section-2add}
In this section we describe iterative EKI as implemented in this
paper. We outline it first for inverse problems as parameterized 
in equations \eqref{inverse} and then discuss generalizations to
centered hierarchical inversion and the non-centered hierarchical
inversion \eqref{hier_inv2}. 
We briefly mention the continuous time limits of the methods as
these provide insight into how EKI works, and the effect of re-parameterizing;
the study of continuous limits from EKI was introduced in \cite{Claudia}
and further details concerning their application to the problems considered
here may be found in \cite{NKC2}.

\subsection{Formulation for \eqref{inverse}}
\label{sec:form1}

The form of iterative EKI that we use is that employed
in \cite{reg}. When applied to the inverse problem \eqref{inverse} it takes the following
form, in which the subscript $n$ denotes the iteration step, and the superscript
$(j)$ the ensemble member: 
\begin{equation}
\label{eq:abasic}
u^{(j)}_{n+1} = u^{(j)}_n + C^{uw}_n(C^{ww}_n + \Upsilon_n\Gamma)^{-1}(y - \mathcal{G}(u^{(j)}_n)).
\end{equation}
The empirical covariances $C^{uw}_n, C^{ww}_n$ are given by 
\begin{align}
C^{uw}_n &= \frac{1}{J-1}\sum_{j=1}^{J}(u^{(j)}_n - \bar{u}_n) \otimes (\mathcal{G}(u^{(j)}_n) - \bar{\mathcal{G}}_n)  \\
C^{ww}_n &= \frac{1}{J-1}\sum_{j=1}^{J} (\mathcal{G}(u^{(j)}_n) - \bar{\mathcal{G}}_n)  \otimes (\mathcal{G}(u^{(j)}_n) - \bar{\mathcal{G}}_n).
\end{align}
Here $\bar{u}_n$ denotes the average of $u^{(j)}_n$ over all ensemble
members and $\bar{\mathcal{G}}_n$ denotes the average of $\mathcal{G}(u^{(j)}_n)$ over all ensemble
members. The parameter $\Upsilon_n$ is chosen to ensure that
\begin{equation}
\label{princp}
\|y-\bar{\mathcal{G}}_n\|_{\Gamma} \leq \zeta \eta,
\end{equation}
a form of discrepancy principle which avoids over-fitting.

\begin{algorithm}
\caption{Hierarchical Iterative Kalman Method (Centered Version).}
\label{algorithm-2}
\begin{algorithmic}[h!] \\
\State Let $\{u_{0}^{(j)}, \theta^{(j)}_{0} \}^{J}_{j=1} \subset \mathcal{X}$ be the initial ensemble with $J$ elements. \\
Further let $\rho \in (0,1)$ with $\zeta > \frac{1}{\rho}$ and $\theta = (\alpha,\tau).$ \\\\ Generate $\{u_{0}^{(j)}, \theta^{(j)}_{0} \}$ i.i.d. from the 
prior $\bbP(u,\theta),$ with synthetic data \\ $y^{(j)}_{n+1}=y + \eta^{(j)}_{n+1}, \ \ \eta^{(j)} \sim N(0,\Gamma)$ i.i.d. \\\\
Then for $n=1,\ldots$
\State \item[1.] \textbf{Prediction step:} Evaluate the forward map $w^{(j)}_{n} = \mathcal{G}(u^{(j)}_{n})$, \\\\
and define $\bar{w}_{n} = \frac{1}{J}\sum^{J}_{j=1}w^{(j)}_{n}$.
\State \item[2.]  \textbf{Discrepancy principle:} If $\| \Gamma^{-1/2}(y-\bar{w}_n)\|_{\mathcal{Y}} \leq \zeta \eta$, {stop!} \\\\ Output $\bar{u}_n = \frac{1}{J}\sum^{J}_{j=1}u^{(j)}_n$ and  $\bar{\theta}_n = \frac{1}{J}\sum^{J}_{j=1}\theta^{(j)}_n$. 
\State \item[3.] \textbf{Analysis step:} Define sample covariances:
\\\\
$ C^{uw}_{n} = \frac{1}{J - 1} \sum^{J}_{j=1}(u^{(j)} - \bar{u}) \otimes (\mathcal{G}(u^{(j)}) - \bar{\mathcal{G}})$,  
\\ 
$ C^{\theta w}_{n} = \frac{1}{J - 1} \sum^{J}_{j=1}(\theta^{(j)} - \bar{\theta}) \otimes (\mathcal{G}(u^{(j)}) - \bar{\mathcal{G}})$,
\\
$ C^{ww}_{n} = \frac{1}{J - 1} \sum^{J}_{j=1}(\mathcal{G}(u^{(j)}) - \bar{\mathcal{G}})  \otimes  (\mathcal{G}(u^{(j)}) - \bar{\mathcal{G}})$.
\\\\
Update each ensemble member as follows \\\\ $u^{(j)}_{n+1} = u^{(j)}_{n} + C^{uw}_n (C^{ww}_n +  \Upsilon_n \Gamma)^{-1}(y^{(j)}_{n+1} - \mathcal{G}(u^{(j)}_n)),$ \\
 $\theta^{(j)}_{n+1} = \theta^{(j)}_{n} + C^{\theta w}_n (C^{ww}_n + \Upsilon_n \Gamma)^{-1}(y^{(j)}_{n+1} - \mathcal{G}(u^{(j)}_n))
$,
\\\\  where $\Upsilon_{n}$ is chosen as $\Upsilon_{n}^{i+1} = 2^{i}\Upsilon ^{0}_{n}$, \\\\
 where $\Upsilon ^{0}_{n}$ is an initial guess. We then define $\Upsilon_n \equiv \Upsilon^N_n$ where $N$ is the first integer such that \\\\
 $\rho \| \Gamma^{-1/2}(y^{(j)}-\bar{w}_n)\|_{\mathcal{Y}} \leq \Upsilon ^{N}_n\| \Gamma^{1/2}(C^{ww}_{n}+\Upsilon ^{N}_n\Gamma)^{-1}(y^{(j)}-\bar{w}_n)\|_{\mathcal{Y}}$. \\\\
\end{algorithmic}
\label{a:1}
\end{algorithm} 

If we define
$$d_n^{(j,m)}=\langle (C^{ww}_n + \Upsilon_n\Gamma)^{-1}(\mathcal{G}(u^{(j)}_n)-y), \mathcal{G}(u^{(m)}_n) - \bar{\mathcal{G}}_n\rangle,$$
then we see that
\begin{equation}
\label{eq:u1}
u^{(j)}_{n+1} = u^{(j)}_n-\frac{1}{J-1}\sum_{m =1}^{J} d_n^{(j,m)} u^{(m)}_n,
\end{equation}
and it is apparent that the algorithm will preserve the linear span of the initial
ensemble $\{u^{(j)}_0\}_{j=1}^{J}.$ We describe the details of how $\Upsilon_n$ is chosen
in the next subsection where we display the algorithm in full for a generalization of the setting of \eqref{inverse} to the hierarchical setting.

To write down the continuous-time limit of the EKI we consider
the setting in which $\Upsilon_n^{-1} \equiv (J-1)h$ and view $u_n^{(j)}$ as
approximating a function $u^{(j)}(t)$ at time $t=nh.$ 
If we define
$$d^{(j,m)}=\langle \Gamma^{-1}(\mathcal{G}(u^{(j)})-y), \mathcal{G}(u^{(m)}) - \bar{\mathcal{G}}\rangle,$$
with obvious definition of $\bar{\mathcal{G}}$,
then we obtain the continuous-time limit 
\begin{equation}
\label{eq:limit1}
\dot{u}^{(j)}=-\sum_{m=1}^J d^{(j,m)}u^{(m)},
\end{equation}
{where $\dot{u}^{(j)}$ denotes the standard time derivative of $u^{(j)}$ viewed as solving an ordinary differential equation; because
the algorithm preserves the linear span of the initial ensemble \cite{Claudia}
the dynamics take place in a finite dimensional space, even if ${\mathcal{X}}$
is inifnite dimensional}. Note that $d^{(j,m)}$ depends on $\{u^{(k)}\}_{k=1}^{J}$ and
that the dynamical system couples the ensemble members. 

\subsection{Generalization for Centered Hierarchical Inversion}

{
Algorithm \ref{a:1} shows the generalization of \eqref{eq:abasic} to the setting
of centered hierarchical inversion. The ensemble is now over both $u$ and $\theta$ and cross
covariances from the observational space to both the $u$ and $\theta$ spaces are
required. 
Algorithm \ref{a:1} also spells out in detail how the parameter $\Upsilon_n$
is chosen.  We now define
$$d_n^{(j,m)}=\langle (C^{ww}_n + \Upsilon_n\Gamma)^{-1}(\mathcal{G}(u^{(j)}_n)-y), \mathcal{G}(u^{(m)}_n) - \bar{\mathcal{G}}_n\rangle,$$
and we see that
\begin{equation}
\label{eq:u1_c}
u^{(j)}_{n+1} = u^{(j)}_n-\frac{1}{J-1}\sum_{m =1}^{J} d_n^{(j,m)} u^{(m)}_n,
\end{equation}
and
\begin{equation}
\label{eq:u2_c}
\theta^{(j)}_{n+1} = \theta^{(j)}_n-\frac{1}{J-1}\sum_{m =1}^{J} d_n^{(j,m)} \theta^{(m)}_n.
\end{equation}
It is apparent that, once again, the algorithm will preserve the linear span of the initial
ensemble $\{u^{(j)}_0,\theta^{(j)}_0\}_{j=1}^{J}.$ 
Furthermore we note that for the centered hierarchical method, since $\mathcal{G}$ does
not depend on $\theta$, the algorithm projected onto the $u$ coordinate
is identical to that
in the preceding subsection, with the only difference being that the initial
span of $\{u^{(j)}_0\}_{j=1}^{J}$ is constructed over a diverse set of $\theta$,
reflecting the dependency structure in $\bbP(u,\theta);$
for hierarchical priors as in subsubsections \ref{sssec:h}
and \ref{subsubsection-2.2.3}, 
the dependency structure is typically 
of the form $\bbP(u,\theta)=\bbP(u|\theta)\bbP(\theta).$}
\footnote{{The  centered hierarchical method, where $\mathcal{G}$ does
not depend on $\theta$, is presented for
simplicity. The details of this algorithm are readily transferred 
to include $\theta$ dependence in the forward mapping $\mathcal{G}$
as required by the version of the Bayesian level set method advocated 
in \cite{MD2}.}}

{
If we again define
$$d^{(j,m)}=\langle \Gamma^{-1}(\mathcal{G}(u^{(j)})-y), \mathcal{G}(u^{(m)}) - \bar{\mathcal{G}}\rangle,$$
with obvious definition of $\bar{\mathcal{G}}$,
then we obtain the continuous-time limit 
\begin{align}
\label{eq:limit2}
\dot{u}^{(j)}=-\sum_{m=1}^J d^{(j,m)}u^{(m)}, \\
\dot{\theta}^{(j)}=-\sum_{m=1}^J d^{(j,m)}\theta^{(m)}.
\end{align}
Since $d^{(j,m)}$ depends only on $\{u^{(k)}\}_{k=1}^{J}$ and
not on $\{\theta^{(k)}\}_{k=1}^{J}$
for the centered hierarchical method
the continuous time limit for $u$ is identical 
to that in the preceding subsection, with the only difference being the creation
of the initial ensemble using variable $\theta.$ This severely limits
the capability of the hierarchical method in the centered case, and is motivation
for the non-centered approach that we now describe.  
}

\subsection{Generalization for Non-Centered Hierarchical Parameterization}

One of the reasons for using hierarchical parameterizations is that good choices of
parameters such as length scale and regularity of $u$ are not known a priori; this suggests
that they might be learnt from the data. In this context the preservation of the linear 
span of the initial ensemble is problematic if the algorithm is formulated in terms
of $(u,\theta).$ This is because, even though the length-scale (for example) may update
as the algorithm progresses, the output for $u$ remains in the linear span of
the initial set of $u$, which likely does not contain a good estimate of the true length
scale. Instead one can work with the variables $(\xi,\theta)$, where $\xi$ {is the forcing function in}
a stochastic PDE, as explained in subsubsections \ref{sssec:h}
and \ref{subsubsection-2.2.3}. 
Working with $(u,\theta)$ and with $(\xi,\theta)$ 
are referred to as the \textit{centered} parameterization and the 
\textit{non-centered} parameterization, respectively. The pros and cons of each
method is discussed in the context of Bayesian inversion in \cite{centre1,centre2},
where the terminology is also introduced. The provenance of the terminology has
no direct relevance in our context, but we retain it to make the link with the
existing literature.

{
The algorithm for updating $(\xi,\theta)$ is identical to that shown in
subsection \ref{sec:form1} with the identifications $u \mapsto (\xi,\theta)$ 
and $\mathcal{G} \mapsto \mathcal{G} \circ T.$} Note that even though, in the
centered case, $\mathcal{G}$ depends only on $u$, the mapping 
$\mathcal{G} \circ T$ will depend on both
$\xi$ and $\theta$. Hence these variables are coupled through the
iteration. Indeed if we now define
$$d_n^{(j,m)}=\langle (C^{ww}_n + \Upsilon_n\Gamma)^{-1}(\mathcal{G}\circ T(\xi^{(j)}_n,\theta^{(j)}_n)-y), \mathcal{G}\circ T(\xi^{(m)}_n,\theta^{(m)}_n) - \overline{\mathcal{G}\circ T}_n\rangle$$
then we see that
\begin{equation}
\label{eq:u1_nc}
\xi^{(j)}_{n+1} = \xi^{(j)}_n-\frac{1}{J-1}\sum_{m =1}^{J} d_n^{(j,m)} \xi^{(m)}_n,
\end{equation}
and
\begin{equation}
\label{eq:u2_nc}
\theta^{(j)}_{n+1} = \theta^{(j)}_n-\frac{1}{J-1}\sum_{m =1}^{J} d_n^{(j,m)} \theta^{(m)}_n.
\end{equation}
Although the algorithm will preserve the linear span of the initial
ensemble $\{\xi^{(j)}_0,\theta^{(j)}_0\}_{j=1}^{J}$ 
the variable of interest $u^{(j)}_n=T(\xi^{(j)}_n,\theta^{(j)}_n)$ is not
in the linear span of $u^{(j)}_0$, in general. 
This confers a significant
advantage on the non-centered parameterization in comparison with
the centered approach.

As in the previous subsections we describe a continuous-time limit, now
for the non-centered hierarchical approach. We define 
$$d^{(j,m)}=\langle (\Gamma^{-1}(\mathcal{G}\circ T(\xi^{(j)},\theta^{(j)})-y), \mathcal{G}\circ T(\xi^{(m)},\theta^{(m)}) - \overline{\mathcal{G}\circ T}\rangle,$$
again with the obvious definition of $\overline{\mathcal{G}\circ T}.$
In the same setting adopted in the previous two subsections we obtain the 
limiting equations
\begin{align*}
\label{eq:limit2}
\dot{\xi}^{(j)}=-\sum_{m=1}^J d^{(j,m)}\xi^{(m)}, \\
\dot{\theta}^{(j)}=-\sum_{m=1}^J d^{(j,m)}\theta^{(m)}.
\end{align*}
Now $d^{(j,m)}$ depends on both $\{\xi^{(k)}\}_{k=1}^{J}$ and
on $\{\theta^{(k)}\}_{k=1}^{J}$
so that the dynamical system not only couples the ensemble members
but in general can couple the dynamics for $\xi$ and for $\theta.$
This is another way to understand the significant
advantage of the non-centered parameterization in comparison with
the centered approach.

\section{Model Problems}\label{section_model}

In order to demonstrate the benefits of the parameterizations that
we introduced here we employ a number of models on which we will base 
our numerical experiments. This section will be dedicated to describing 
the various PDEs that will be used. {We will describe the forward 
problem, together with a basic version of the inverse problem, for each model, 
relevant in the non-hierarchical case. We note that the ideas such as the
level set method, and hierarchical formulations from 
subsubsections \ref{subsubsection-2.1.2} and \ref{subsubsection-2.1.3}, can 
be used to reformulate the inverse problems, and we will use these
reformulations in section \ref{section-4}.}

\subsection{Model Problem 1}\label{subsection-3.4}

Our first test model is from electrical impedance tomography (EIT). 
This imaging method is used to learn about interior properties
of a medium by injecting current, and measuring voltages, on the
boundary \cite{EIT,MD}.  We will use the complete electrode model (CEM)
introduced in \cite{EIT2}. The forward model is as follows:
{given a domain $D = \mathcal{B}(0,1)^2$ and a set of electrodes $\{e_l\}_{1=1}^{m_e}$ on the boundary $\partial D$} with contact impedance $\{z_l\}_{l=1}^{m_e}$, and interior conductivity $\kappa$, the CEM aims to solve for the potential $\nu$ inside the domain $D$ and the voltages 
$\{V_l\}_{1=1}^{m_e}$ on the boundary. The governing equations are 
\begin{subequations}
\label{CEM}
\begin{alignat}{4}
\nabla \cdot (\kappa \nabla \nu) &= 0 , \   \ \in \ D \\
\nu + z_l\kappa \nabla\nu \cdot n &= V_l , \   \in \ e_l,  \ \ l=1,\ldots,m_{e}  \\
\nabla\nu \cdot n &=0 , \  \ \in \ \partial D \textrm{\textbackslash} \cup^{m_{e}}_{l=1}e_{l}  \\
\int\kappa \nabla \nu \cdot n \ ds &= I_l , \   \in \ e_l,  \ \ l=1,\ldots,m_{e},
\end{alignat}
\end{subequations}
with $n$ denoting the outward normal vector on the boundary. The linearity
of the problem implies that the relationship between injected current
and measured voltages can be described through an Ohm's Law
of the form
\begin{equation}
V = R(\kappa) \times I.
\end{equation}
In our experiments $D$ will be a two-dimensional disc of radius $1$.
 The inverse problem may now be stated. We write the unknown conductivity
as $ \kappa = \exp(u)$ and try to infer $u$ from a set of $J$ noisy measurements 
of voltage/current pairs $(V_j,I_j).$ If we define 
$\mathcal{G}_{j}(u) = R(\kappa) \times I_j$ then the inverse problem is to
find $u$ from $y$ given an equation of the form \eqref{inverse}.

\subsection{Model Problem 2}\label{subsection-3.5}
Our second model problem arises in hydrology: 
the single-phase Darcy flow equations. The concrete instance of the
forward problem is as follows: given the  domain $D = [0,6]^2$ and 
real-valued permeability function $\kappa$ defined on $D$,  the forward model 
is to the determine real-valued pressure (or hydraulic head) function
$p$ on $D$ from
\begin{equation}
\label{eq:darcy}
-\nabla\cdot({\kappa}\nabla p) = f, \  \ x \ \in \ D, \\
\end{equation}
with mixed boundary conditions
\begin{equation}
p(x_1,0) = 100, \ \ \ \  \frac{\partial p}{\partial x_1}(6,x_2) = 0, \ \ -\kappa\frac{\partial p}{\partial x_1}(0,x_2) = 500, \ \ \  \frac{\partial p}{\partial x_2} (x_1,6) = 0,
\end{equation}
and the source term $f$ defined as
\[
f(x_1,x_2)=
\begin{cases}
0, & \textrm{if} \ \  0 \leq x_2 \leq 4, \\
137, & \textrm{if} \ \  4 \leq x_2 \leq 5, \\
274,  & \textrm{if} \ \  5 \leq x_2 \leq 6.
\end{cases}
\]
The inverse problem concerned with (\ref{eq:darcy}) is as follows:
write $\kappa = \exp(u)$ and determine $u$ from $J$ linear functionals 
of the pressure $\mathcal{G}_{j}(u) = l_{j}(p)$. 
This may thus be cast in the form \eqref{inverse}.
We take the linear functionals as mollified pointwise observations on a regular
grid.
This specific set-up of the PDE model is that tested by Hanke \cite{Hanke} in 
his consideration of the regularized Levenberg-Marquardt algorithm. More information on this setting can be found by Carera et al. \cite{BC}.

\subsection{Model Problem 3}\label{subsection-3.6}
Our final model is a simple linear inverse problem
which we can describe directly. The aim is to reconstruct a function 
$u$ from noisy observation of $J$ linear functionals $\mathcal{G}_{j}(u) = l_{j}(p)$, $j=1,\cdots, J$, where $p$ solves
the equation
\begin{subequations}
\label{1D}
\begin{align}
\frac{d^2p}{dx^2} + {p} &= u,\ \in D, \\
p&=0 , \ \in \partial D.
\end{align}
\end{subequations}
This may also be cast in the form \eqref{inverse}. 
We use equally spaced pointwise evaluations as our linear functionals. {We will assume our domain is chosen such that $D=[0,10]$.}

\section{Numerical Examples} \label{section-4}

To assess the performance of each parameterization we present a range
of numerical experiments on each of the three model problems
described in the previous section. Our experiments will be presented
in a consistent fashion, between the different models and the
different algorithms.  Each model problem will be tested using each
of the non-hierarchical and hierarchical approaches, although we will not
use the centered approach for Model Problem 3. 
Within each of these approaches we will show the progression 
of the inverse solver from the first to the last iteration. This will include 
five images ordered by iteration number, with the first figure displaying
the first iteration and the last displaying the final iteration. These figures 
will be accompanied with figures demonstrating the learning of the
hyperparameters, as the iteration progresses, for Model Problems 1 and 2,
but not for Model Problem 3 (where the hyperparameter is a field). 

In order to illustrate the effect of the initial ensemble we will show output
of the non-centered approach for ten different initializations, for each
model problem.  We will plot the final iteration reconstruction arising 
from four of those initializations. {We observe the variation across
the initializations through the relative errors in the unknown field $u_{{\tiny{\rm EKI}}}$, with {respect to} the truth $u^{\dagger}$, and in the data misfit, as
the iteration progresses:} \\
\begin{align*}
{\frac{\| u_{{\tiny{\rm EKI}}} - u^{\dagger}\|_{L^2(D)}}{ \| u^{\dagger} \|_{L^2(D)}}}, \\
{\| y - \bar{\mathcal{G}}(u_{\rm EKI})\|_{\Gamma}}.
\end{align*}

\subsection{Level Set Parameterization}\label{subsection-4.1}
Level set methods are a computationally effective way to represent
piecewise constant functions, and there has been considerable
development and application to inverse problems \cite{Burger1, Burger2, reg},
starting from the paper \cite{Santosa1996}, in which interfaces
are part of the unknown.
We apply level set techniques, combined with hierarchical parameter estimation,
in the context of ensemble inversion; we are
motivated by the recent Bayesian level set method developed by Lu et al. 
in \cite{Lu}, and its hierarchical extensions introduced in \cite{VC1,MD2}.

When applying the level set technique to inverse problems of the form \eqref{inverse}, we modify our forward operator to
\begin{equation}
\label{ls_map}
\mathcal{G} = \mathcal{O} \circ {G} \circ F,
\end{equation}
where $G: X \mapsto Y$ maps the coefficient of the PDE to its solution,
$\mathcal{O}:Y \rightarrow \mathcal{Y}$ is our observational operator, 
and $F:\mathcal{X} \rightarrow X$ is the level set map described by
\begin{equation}
\label{ls_map}
(Fu)(x) \rightarrow \kappa(x) = \sum^{n}_{i=1} \kappa_{i} \mathbbm{1}_{D_i}(x).
\end{equation}
The sets $\{D_i\}_{i=1}^n$ are $n$ disjoint subdomains with union $D$
and whose boundaries define the interfaces. The boundaries are assumed to
be defined through a continuous real-valued function $u$ on $D$ via its
level sets. In order to model the level set method hierarchically we will base our reconstructions on the approaches taken in subsection \ref{subsubsection-2.2.3}.
In general it can be helpful to re-scale the level values as the hierarchical
parameter $\tau$ is learned \cite{MD2}; however if the unknown is binary and the
level set taken as zero, as used in our numerical experiments here,
then this is not a consideration.

We apply level set inversion to Model Problem 1 (EIT) from
subsection \ref{subsection-3.4}. We reconstruct a binary field and
the variable $u$ is, rather than the logarithm of the conductivity $\kappa$,
the level set function defining \eqref{ls_map}: specifically 
the level set formulation is achieved through representing the conductivity
as
\begin{equation}
\kappa(x)=(Fu)(x) = \kappa_{-} \chi_{u \leq 0} + \kappa_{+}\chi_{u >0},
\end{equation}
where $\chi_{A}$ denotes a characteristic function of $A$ with $\kappa_{-}$ and  $\kappa_{+}$ being known positive constants that help define low and high levels of our diffusion coefficient.

\begin{figure}[h]
\centering
 \includegraphics[scale=0.3]{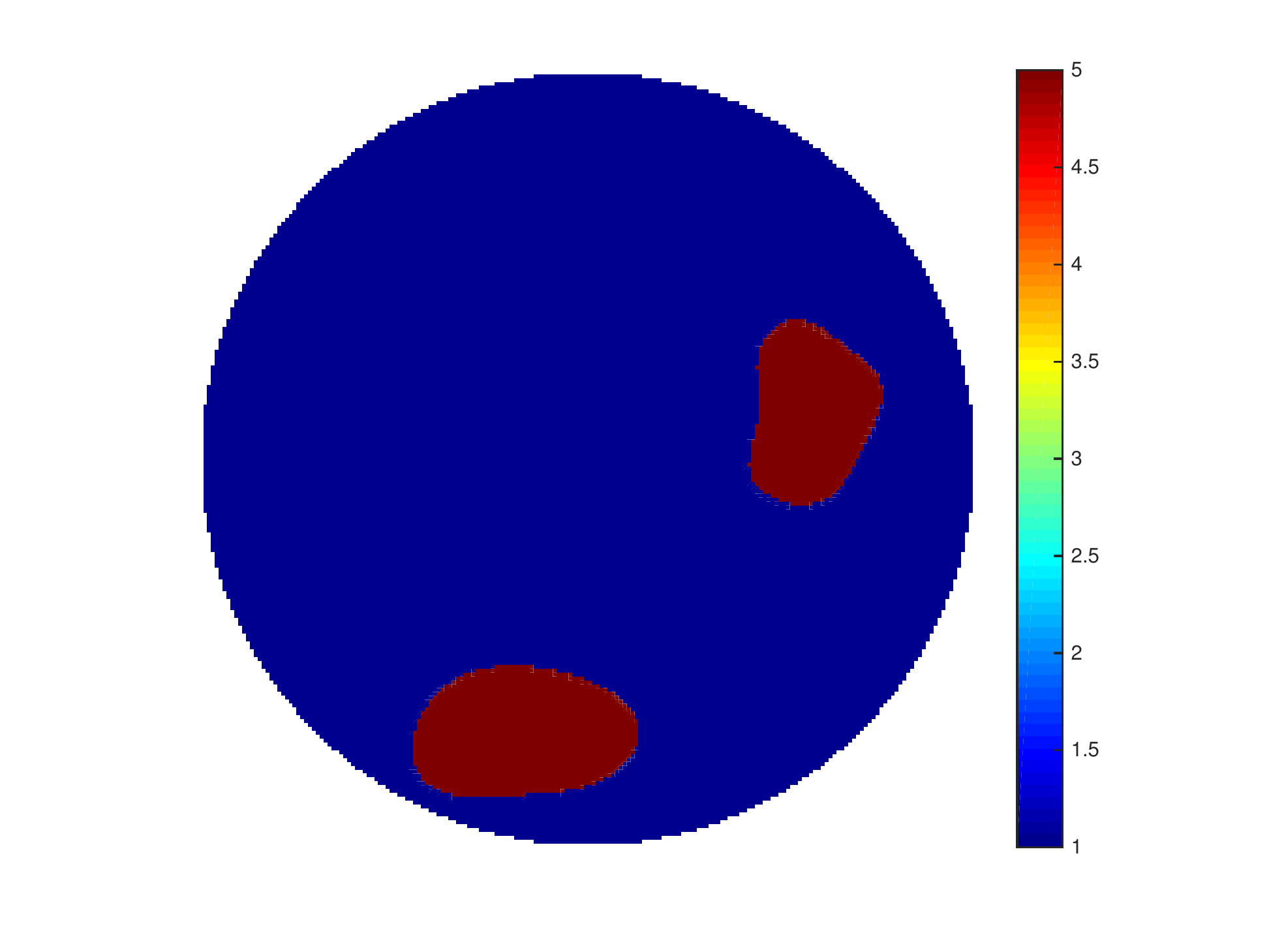}
 \caption{Model problem 1: true log-conductivity.}
 \label{fig:EIT_hier_truth}
\end{figure}

We place $16$ equidistant electrodes on the boundary of the unit disc $D$
in order to define our observations. All experiments are conducted using
the MATLAB package EIDORS \cite{uses}. The contact impedances $\{z_l\}_{l=1}^{m_e}$ are chosen with value $0.05$ and all electrodes chosen subjected to
an input current of $0.1$. This provides a matrix of stimulation patterns $I = \{ I^{(j)}\}_{j=1}^{15} \in \mathbb{R}^{16 \times 15}$ given as

\[
I= 0.1 \times
  \begin{pmatrix}
  
    +1 & 0 & \ldots & 0  \\
    -1 & +1 & \ldots & 0 \\
    0 & -1 & \ddots & 0 \\
     \vdots & \vdots & \ddots & +1 \\
    0 & 0 & 0 & -1 
  \end{pmatrix}.
\]
\\
 For our iterative method we choose $J=200$ ensemble members with regularization parameter $\rho = 0.8$. The covariance of our noise $\eta$ is chosen such that $\Gamma = 10^{-4} \times I$. Our truth for the EIT problem will take the form given in Figure \ref{fig:EIT_hier_truth} where we have high levels of conductivity within the two inclusions. This is constructed by thresholding a
Whittle-Mat\'{e}rn Gaussian random field defined by \eqref{SPDE}, \eqref{cov_prior}; 
true values for the hierarchical parameters used are shown in 
Table \ref{table:level_truth}.

{
\begin{Remark}
We do not display the underlying Gaussian random field $u$ which is 
thresholded to obtain the true conductivity in Figure 5 as this Gaussian 
random field cannot be expected to be reconstructed accurately, in general. 
Furthermore it is important to appreciate that in general  
a true conductivity will not be constructed by such thresholding; the field 
$u$ is simply an algorithmic construct. We do however show $u$, and its 
evolution, in the algorithm, because this information highlights the roles 
of the length scale and regularity parameters.
\end{Remark}
}
When performing inversion we sample
initial ensembles using  the prior distributions shown in  
Table \ref{table:level_dist}. We have set our prior distributions in such a way 
that the true value for each hyperparameter lies within the range specified.

\begin{table}[h!]
\centering
\begin{tabular}{ c c}
\hline
 Hyperparameter &  Value  \\ [0.5ex]
\hline
$\alpha^{\dagger}$ & $3$ \\
$\tau^{\dagger}$ & $10$ \\[0.5ex]
 \hline
\end{tabular}
\caption{Model problem 1. True values for each hyperparameter.}
\label{table:level_truth}
\end{table}

\begin{table}[h!]
\centering
\begin{tabular}{ c c}
\hline
 Hyperparameter &  Prior  \\ [0.5ex]
\hline
$\alpha$ & $\mathcal{U}[1.3,4]$ \\
$\tau$ & $\mathcal{U}[5,30]$ \\[0.5ex]
 \hline
\end{tabular}
\caption{Model problem 1. Prior distribution for each hyperparameter.}
\label{table:level_dist}
\end{table}

\newpage

\begin{figure}[h!]
\centering
\includegraphics[scale=0.8]{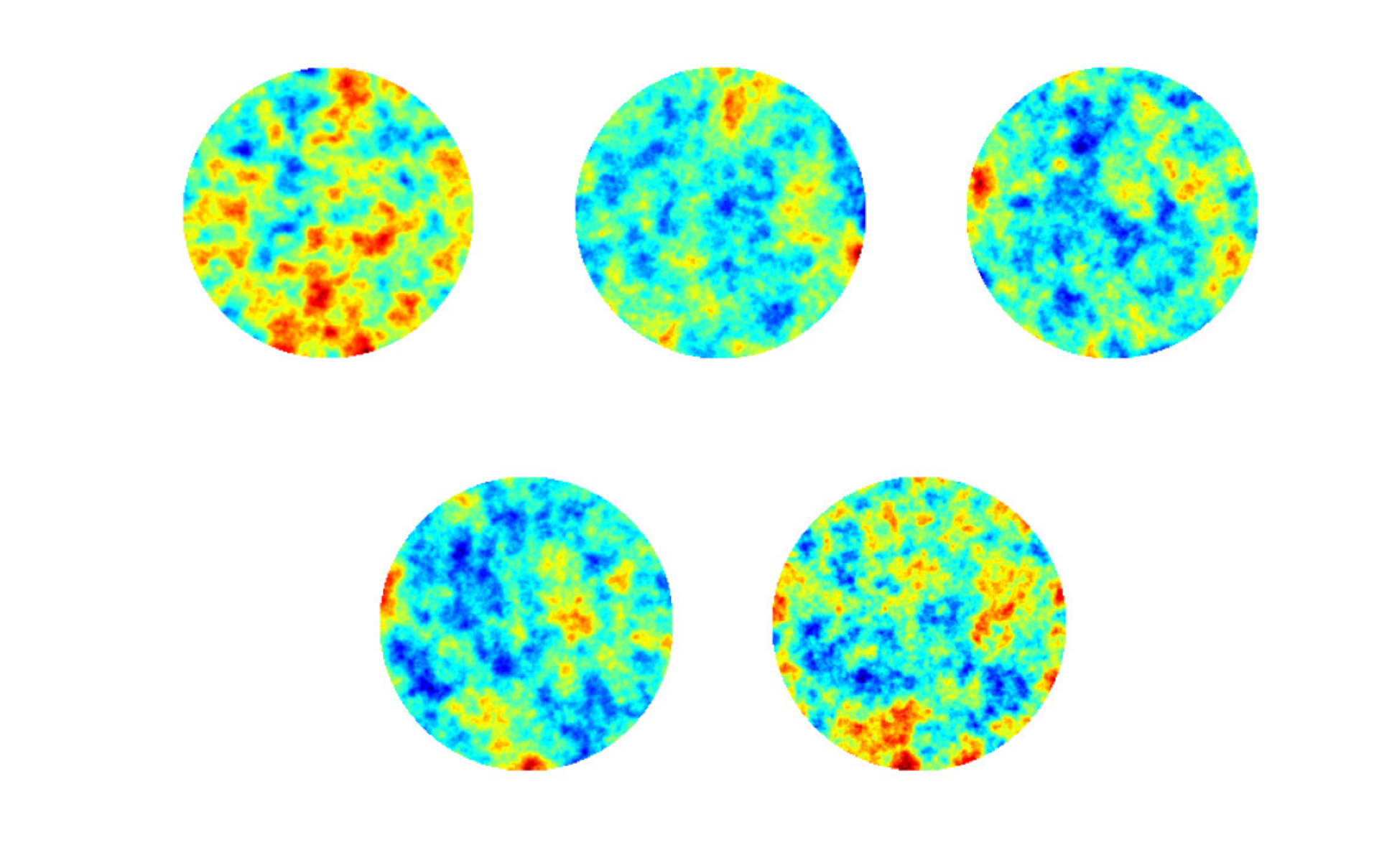}
 \caption{Model problem 1. Progression through iterations of non-hierarchical method.}
 \label{fig:ls_nh_1}
\end{figure}

\begin{figure}[h!]
\centering
 \includegraphics[scale=0.8]{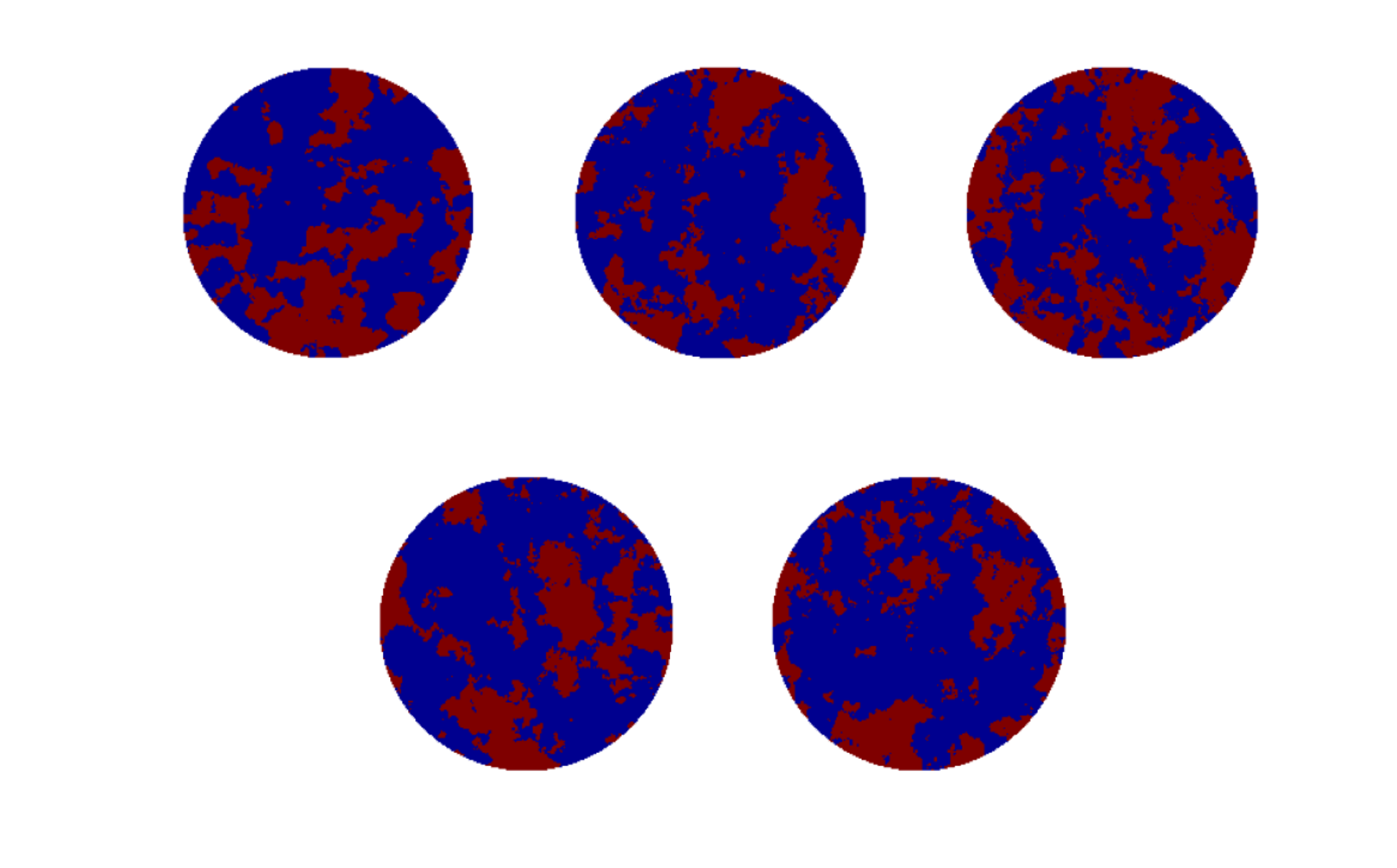}
 \caption{Model problem 1. Progression through iterations of non-hierarchical method with level set.}
 \label{fig:ls_nh_2}
\end{figure}

\begin{figure}[h!]
\centering
\includegraphics[scale=0.8]{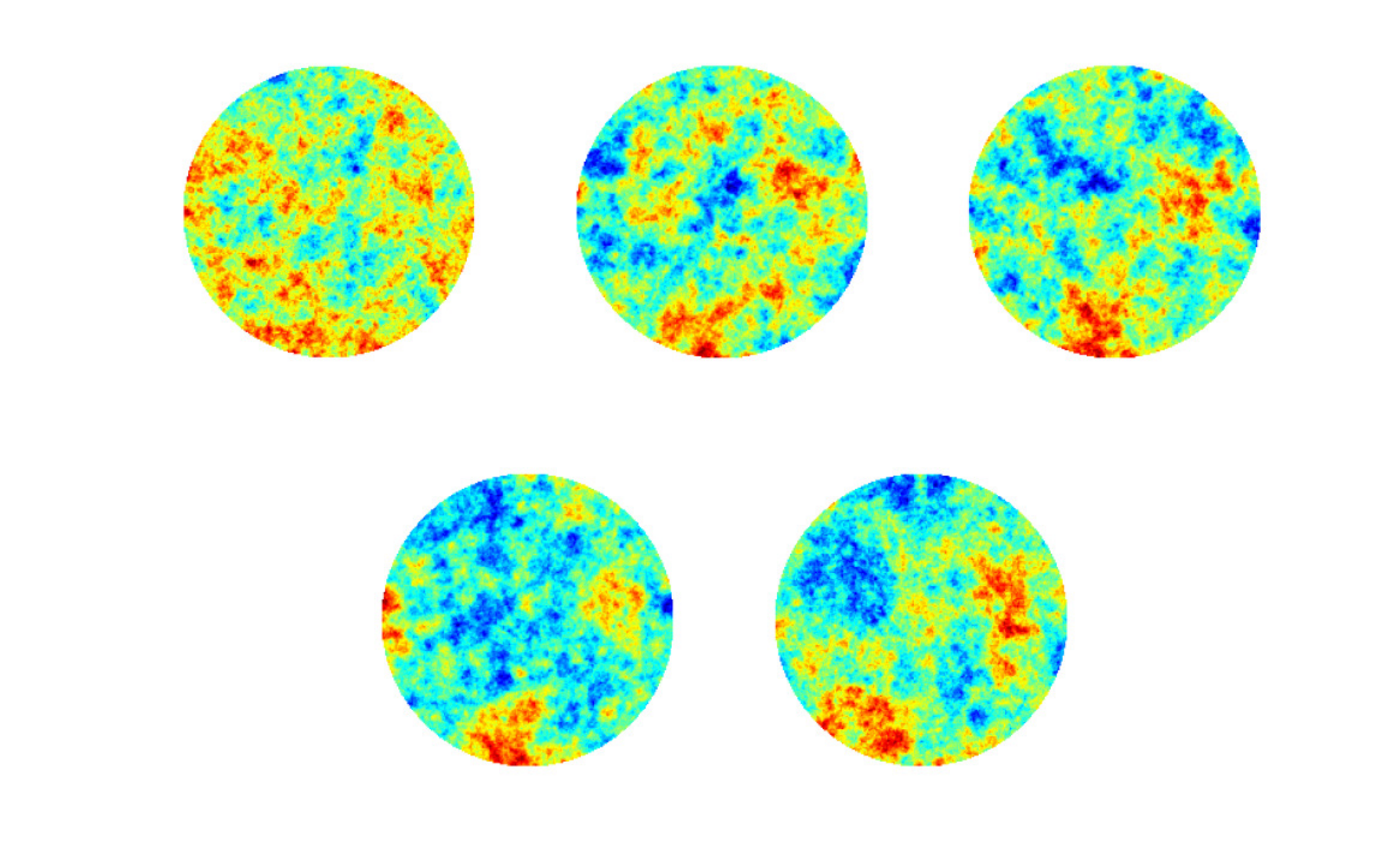}
 \caption{Model problem 1. Progression through iterations of centered hierarchical method.}
 \label{fig:ls_c_1}
\end{figure}

\begin{figure}[h!]
\centering
 \includegraphics[width=\linewidth]{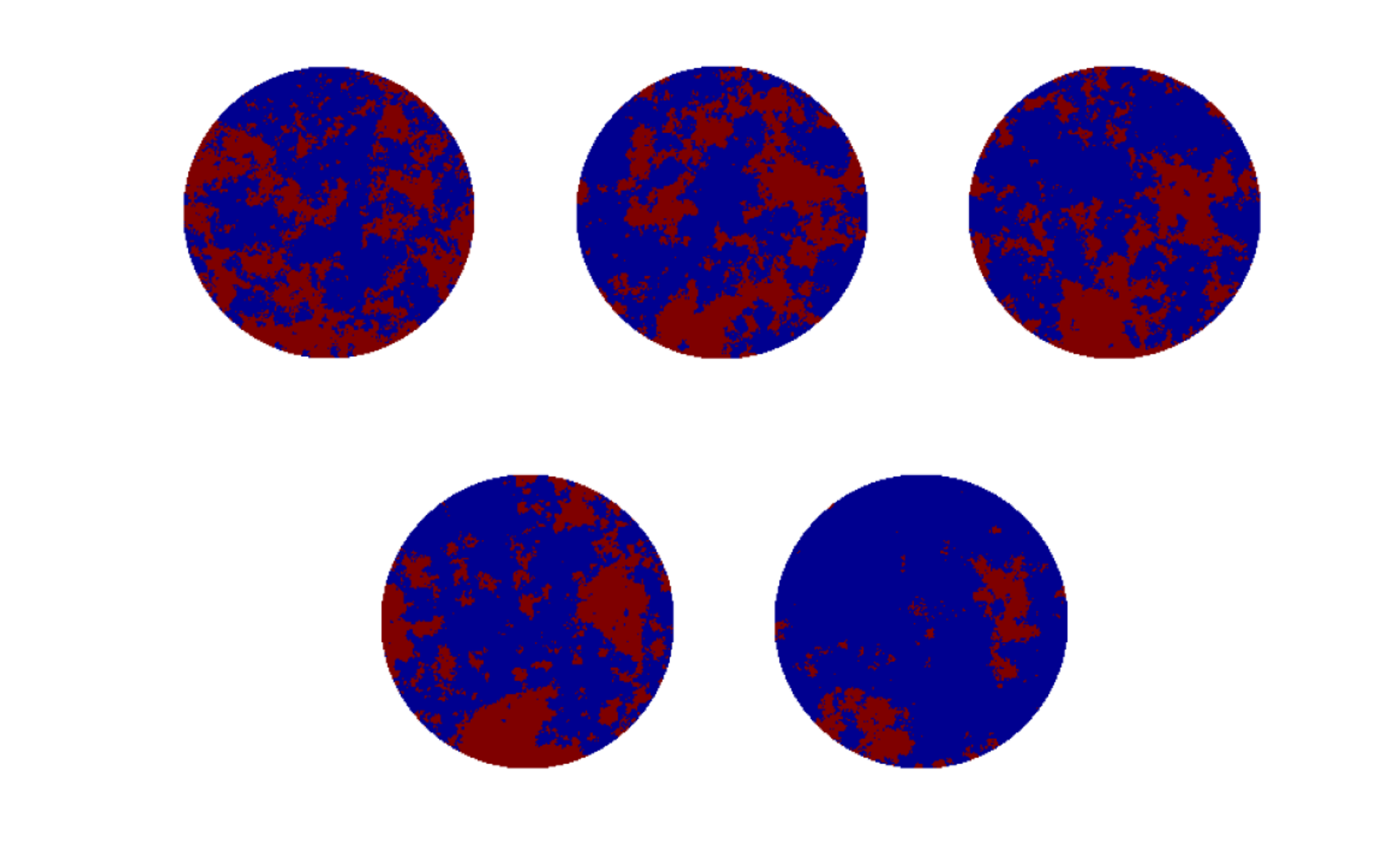}
 \caption{Model problem 1. Progression through iterations of centered hierarchical method with level set.}
 \label{fig:ls_c_2}
\end{figure}

\begin{figure}[h!]
\centering
\includegraphics[width=\linewidth]{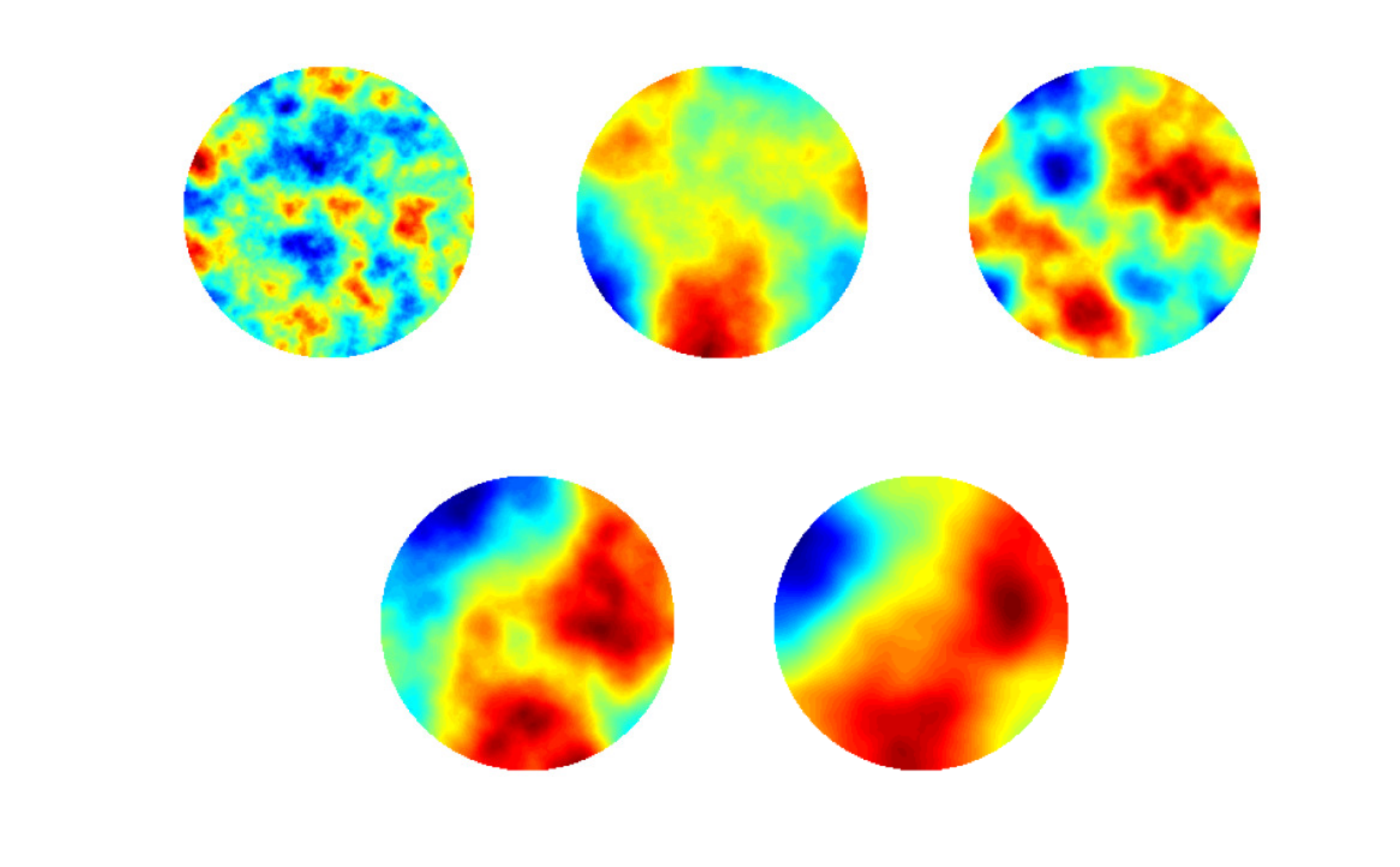}
 \caption{Model problem 1. Progression through iterations of non-centered hierarchical method.}
 \label{fig:ls_nc_1}
\end{figure}

\newpage

\begin{figure}[h!]
\centering
 \includegraphics[width=\linewidth]{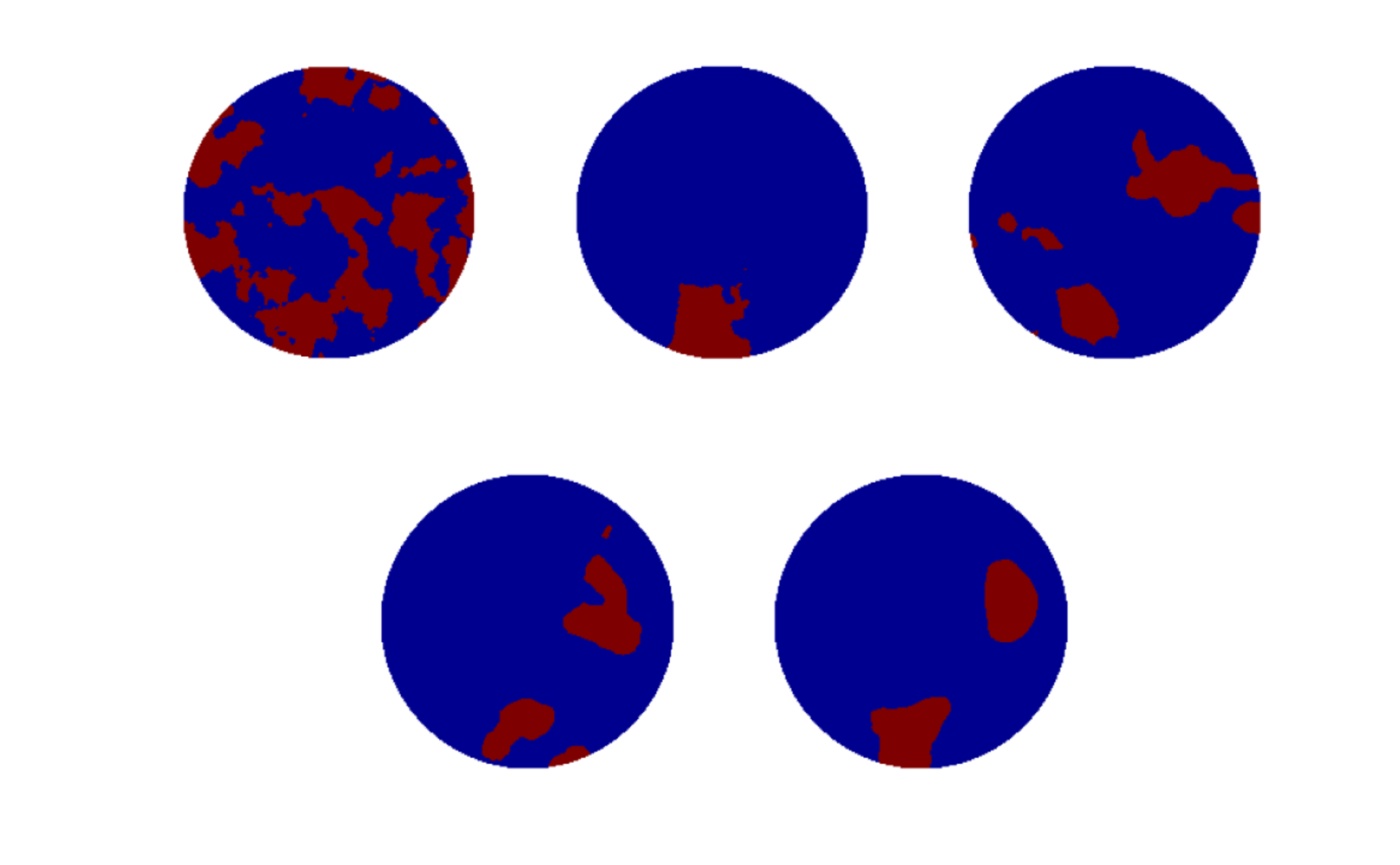}
 \caption{Model problem 1. Progression through iterations of non-centered method with level set.}
 \label{fig:ls_nc_2}
\end{figure}

\begin{figure}[h!]
\centering
 \includegraphics[scale=0.75]{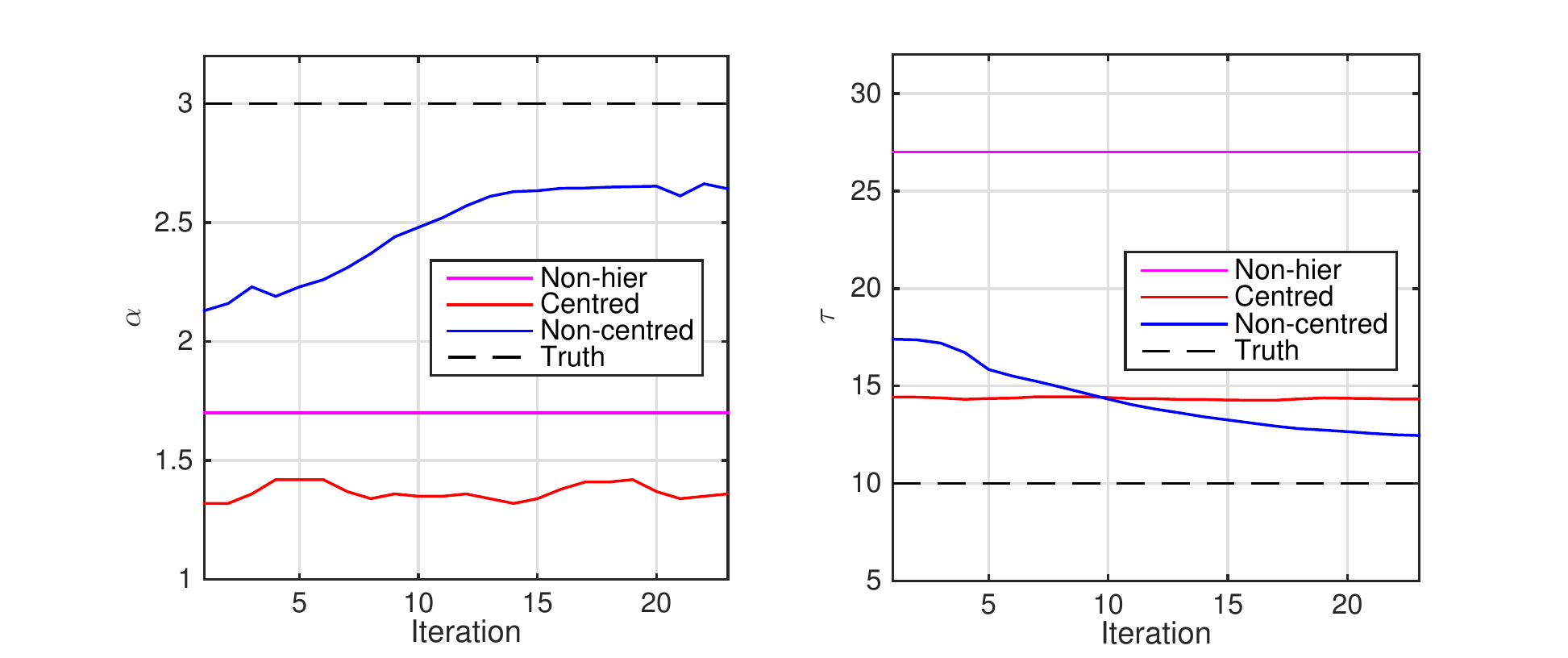}
\caption{Model problem 1. Progression of average value for $\alpha$ and $\tau$.}
 \label{fig:new3}
\end{figure}

\newpage

\begin{figure}[h!]
\centering
 \includegraphics[scale=0.75]{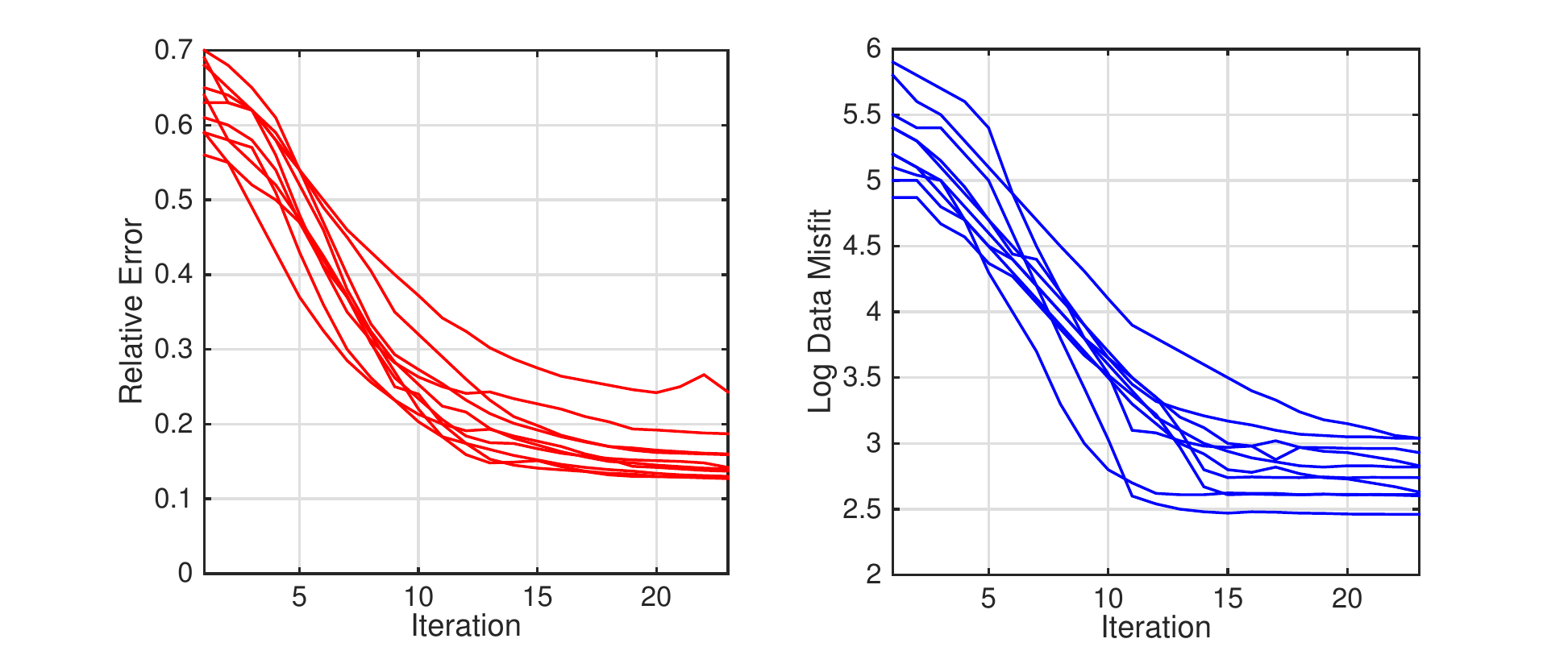}
\caption{Model problem 1. Left: relative error. Right: log-data misfit.}
 \label{fig:levelset_error_misfit}
\end{figure}

\begin{figure}[h!]
\centering
 \includegraphics[scale=0.67]{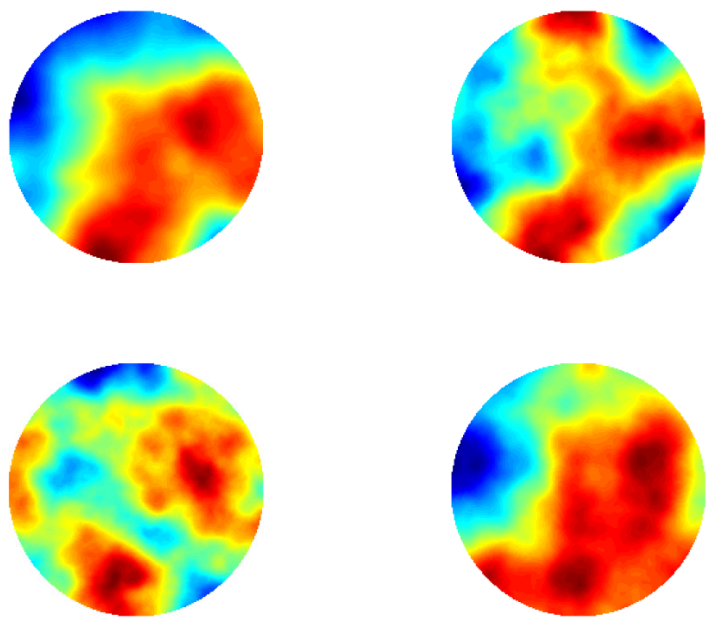}
 \caption{Model problem 1. EKI for the final iteration for the non-centered approach with Whittle-Mat\'{e}rn from four different initializations.}
 \label{fig:final_nc_wm_eit}
\end{figure}

\begin{figure}[h!]
\centering
 \includegraphics[scale=0.67]{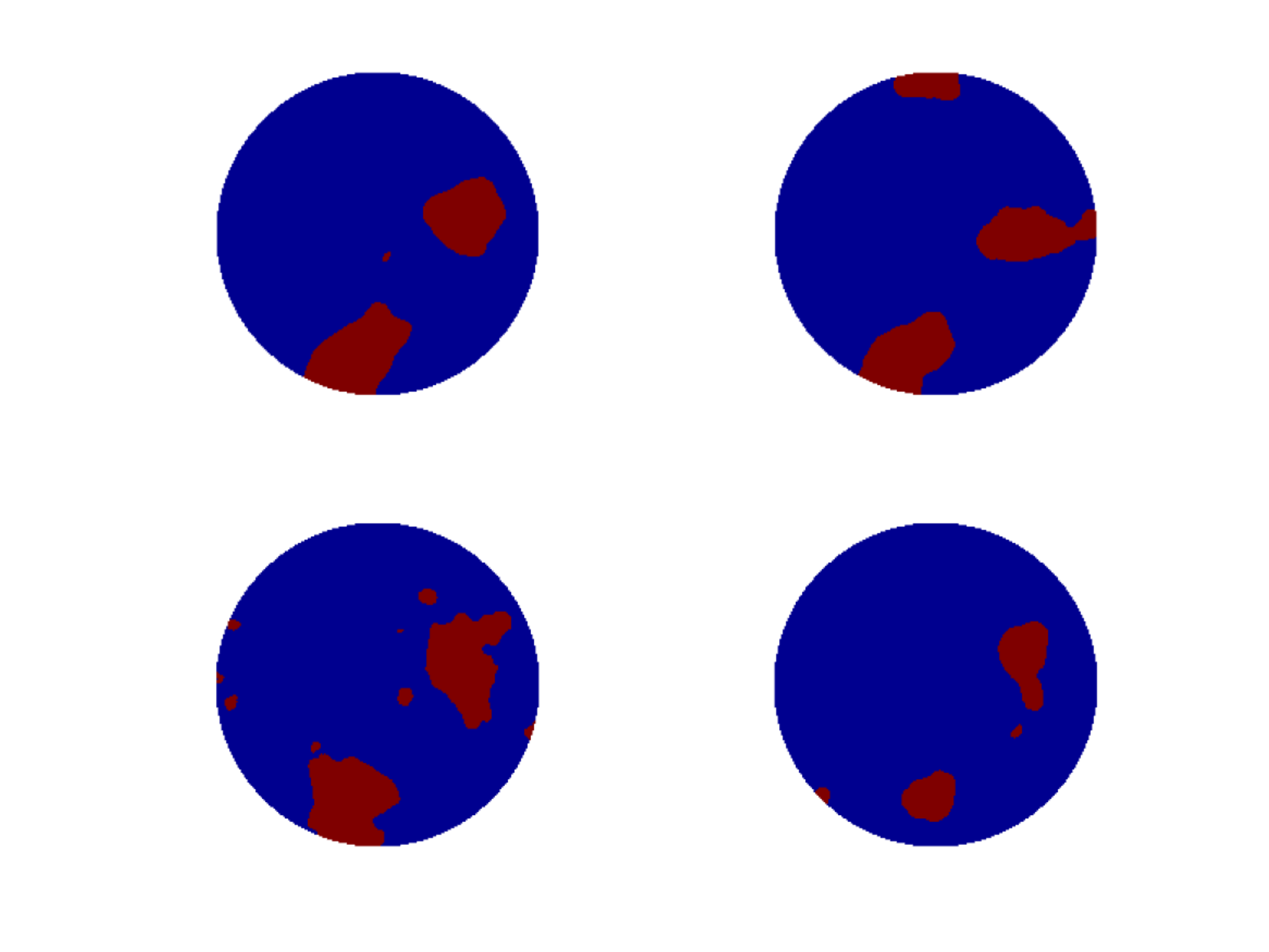}
 \caption{Model problem 1. EKI for the final iteration for the non-centered method with level set from four different initializations. }
 \label{fig:final_nc_ls_eit}
\end{figure}

\newpage

Figures \ref{fig:ls_nh_1} - \ref{fig:ls_nc_2} show the progression of both 
the level-set function $u$, and the permeability $\kappa$, through five 
iterations of the method. Non-hierarchical, centered and non-centered
methods are considered in turn. By comparing the reconstructions with the
true conductivity, these figures clearly demonstrate two facts:
(a) that being hierarchical is necessary to obtain a good reconstruction;
(b) that implementing the hierarchical method using a non-centered 
parameterization has significant benefits when compared to the centered
method. These points are further demonstrated in Figures \ref{fig:new3} 
which shows the learning of the hyperparameters, in comparison with
the truth, for the three methods.

Concentrating solely on the non-centered approach, we run ten different 
initializations of the EKI. 
We display the resulting data-misfit and error, as a function of
iteration, for all ten  in Figure \ref{fig:levelset_error_misfit}. 
We display the last iteration of four of these
ten in Figures \ref{fig:final_nc_wm_eit} -  \ref{fig:final_nc_ls_eit}.

In summary, 
{Figures  \ref{fig:ls_nc_1} -  \ref{fig:new3} clearly show the superiority of the non-centered hierarchical method. For all the initializations shown, the EKI produces conductivities which
concentrate near to the true conductivity and have length scale similar to those
appearing in the truth; see Figures  \ref{fig:levelset_error_misfit} - \ref{fig:final_nc_ls_eit}. The centered hierarchical and non-hierarchical methods fail to do this; see Figures  \ref{fig:ls_nh_1} -  \ref{fig:ls_c_2}. However it is important to note the non-centered method does produce substantial variation in the predicted solution, depending on which initialization is used as shown in Figures   \ref{fig:levelset_error_misfit} - \ref{fig:final_nc_ls_eit}.}



\newpage

\subsection{Geometric Parameterization}\label{subsection-4.2}

In this subsection we employ Model Problem 2 from subsection 
\ref{subsection-3.5}. We consider reconstruction 
of a piecewise continuous channel which is defined through two heterogeneous 
Gaussian random fields, scalar geometric parameters specifying the
geometry and scalar hierarchical parameters characterizing the
length-scale and regularity of the two fields. 

\begin{figure}[h!]
\centering
 \includegraphics[scale=0.25]{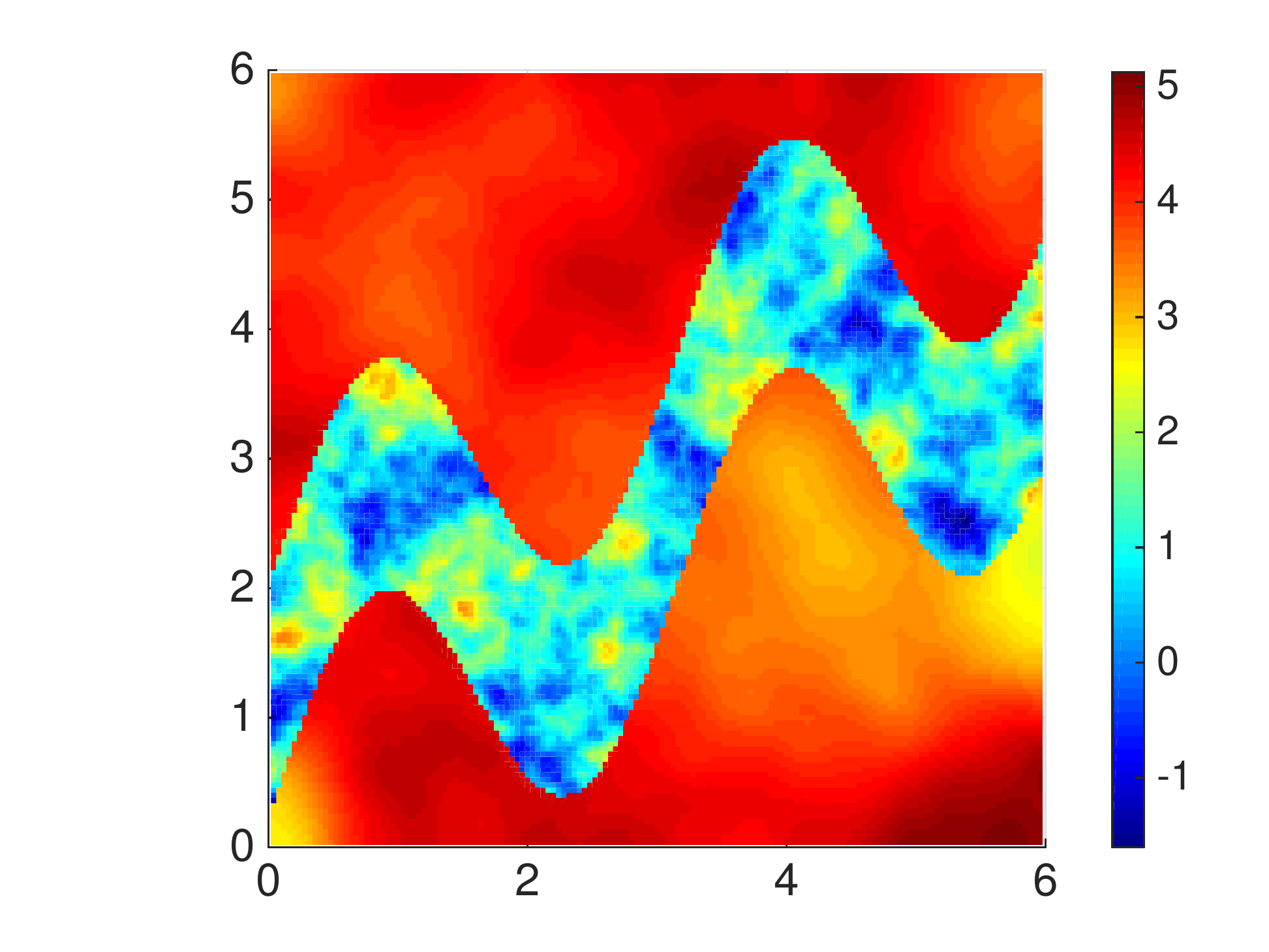}
\caption{Model problem 2. True log-permeability.}
 \label{fig:truth_exp}
\end{figure}

The truth $u^{\dagger}$ is shown in Figure \ref{fig:truth_exp}.
It is drawn from a prior distribution which we now describe; details
may be found in \cite{geo}.
{The channel is described by five parameters: $d_1$ -- amplitude; $d_2$ -- 
frequency; $d_3$ -- angle; $d_4$ -- initial point; and $d_5$ -- width}. 
{We generate two Gaussian random fields $\{\kappa_i\}_{i=1}^2$, both defined
on the whole of the domain $D$ but entering the permeability $\kappa$
only inside and outside (respectively) the channel. The unknown $u$ thus comprises the five scalars
$\{d_i\}_{i=1}^{5}$ and the two fields $\{\kappa_i\}_{i=1}^{2}.$ We do not
explicitly spell out the mapping from $u$ to the coefficient $\kappa$ appearing
in the Darcy flow, but leave this to the reader. The $\{\kappa_i\}_{i=1}^2$ are specified
as log-normal random fields and the underlying Gaussians 
are of Whittle-Mat\'{e}rn type, 
defined by \eqref{SPDE}, \eqref{cov_prior}; different
uniform distributions on $\alpha$ and $\tau$ are used for the two fields
$\{\kappa_i\}_{i=1}^2$ The parameters $\{d_i\}_{i=1}^{5}$ are also given uniform
distributions.} The
entire specification of the prior is given in Table \ref{table:hier_geo_prior}.
For the truth the true hierarchical parameters are provided in Table \ref{table:hyper_dist2}.

In our inversion we employ $64$ mollified pointwise observations $\{ l_{i}(p)\}_{i=1}^{64}$ given by, for some $\sigma >0$, 
\begin{equation}
l_{t}(p) = \int_{D}\frac{1}{2\pi\sigma^2} e^{-\frac{1}{2\sigma^2}(x-x_t)^2} p(x) dx,
\end{equation} 
where the $x_{i}$ are uniformly distributed points on $D$. 
{We discretize the forward model using a second order 
centered finite difference method with mesh spacing $10^{-2}$}. 
For our EKI method we use the same 
values for our parameters as in subsection \ref{subsection-4.1}. 

\begin{table}[ht]
\centering
\begin{tabular}{ c c}
\hline
 Parameter & Prior  \\ [0.5ex]
\hline
    $d_{1}$ &  $\mathcal{U}[0,1]$\\ 
   $d_{2}$ &  $\mathcal{U}[2,13]$ \\
  $d_{3}$ & $\mathcal{U}[0.4,1]$\\
  $d_{4}$ & $\mathcal{U}[0,1]$\\
  $d_{5}$ & $\mathcal{U}[0.1,0.3]$\\    $\kappa_1$ & $N(1, (I- \tau_1^{2} \Delta )^{- \alpha_{1}})$\\
    $\kappa_2$ & $N(4, (I- \tau_2^{2} \Delta )^{- \alpha_{2}})$\\
    $\alpha_1$ & $\mathcal{U}[1.3,3]$\\
    $\tau_{1}$ & $\mathcal{U}[8,30]$\\
    $\alpha_2$ & $\mathcal{U}[1.3,3]$\\
    $\tau_{2}$ & $\mathcal{U}[8,30]$\\  [1ex]
\hline
\end{tabular}
\caption{Model problem 2. Prior associated with channelised flow.}
\label{table:hier_geo_prior}
\end{table}
\bigskip 
\begin{table}[h!]
\centering
\begin{tabular}{ c c c c}
\hline
 Parameter &  Value  \\ [0.5ex]
\hline
$\alpha_1^{\dagger}$ & 2  \\
$\alpha_2^{\dagger}$ & 2.8 \\
$\tau_1^{\dagger}$ &  30\\
$\tau_2^{\dagger}$ & 10 \\[0.5ex]
 \hline
\end{tabular}
\caption{Model problem 2. Parameter selection of the truth.}
\label{table:hyper_dist2}
\end{table}

\newpage

\begin{figure}[h!]
\centering
 \includegraphics[scale=0.75]{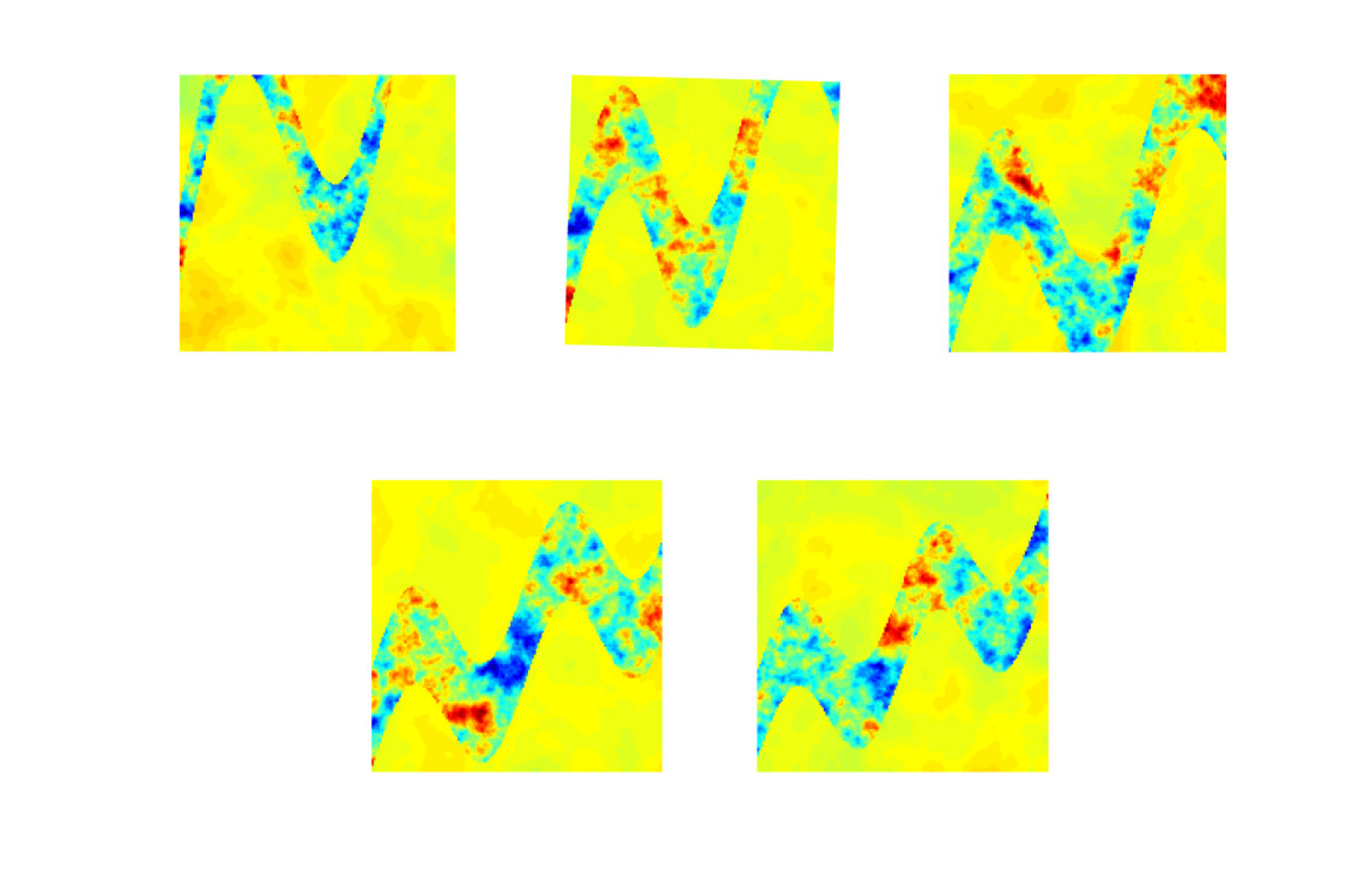}
\caption{Model problem 2. Progression through iterations of non-hierarchical method.}
 \label{fig:non_hier_geo}
\end{figure}

\begin{figure}[h!]
\centering
 \includegraphics[scale=0.75]{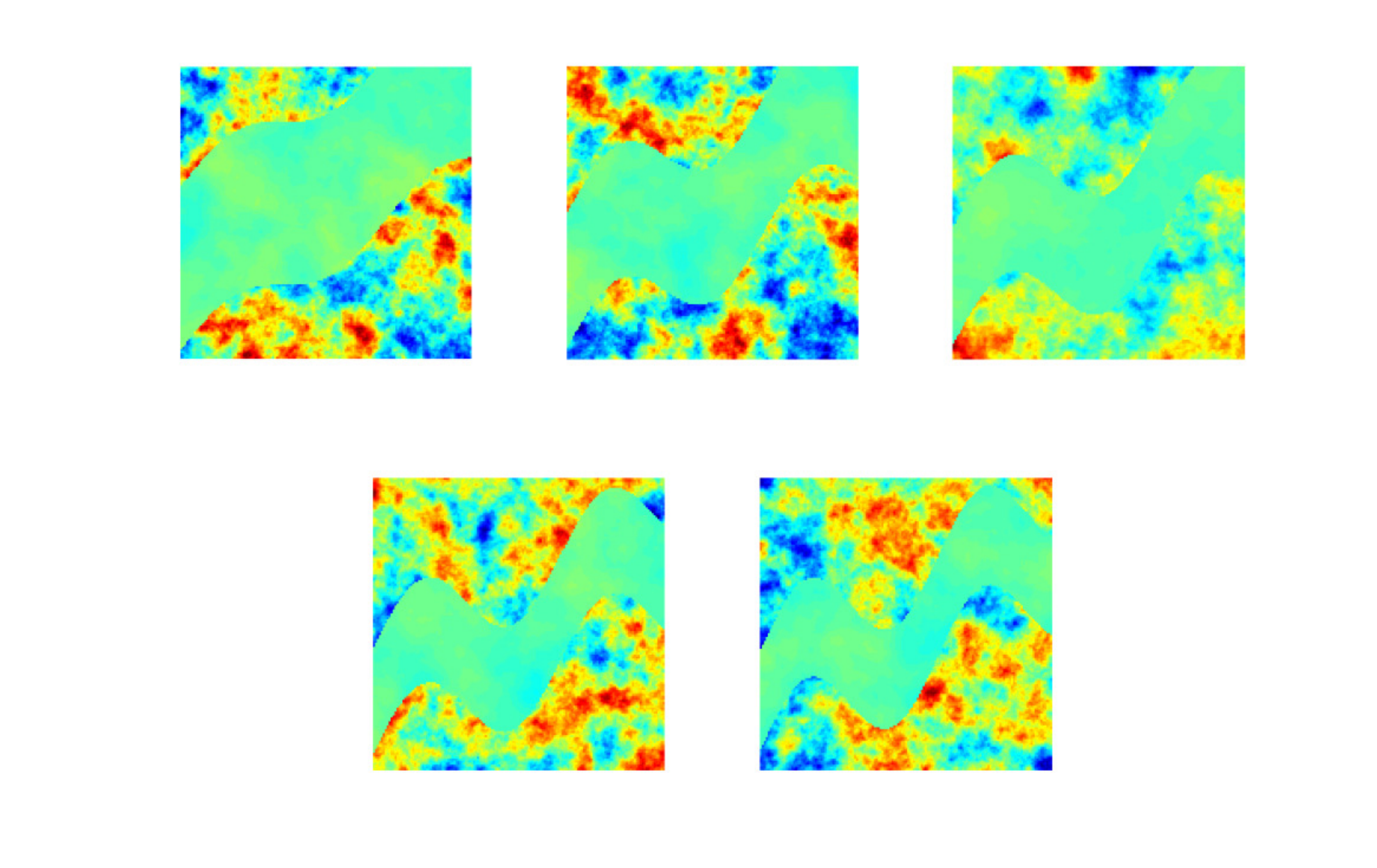}
\caption{Model problem 2. Progression through iterations of centered hierarchical method.}
 \label{fig:hier_c_geo}
\end{figure}

\begin{figure}[h!]
\centering
 \includegraphics[scale=0.75]{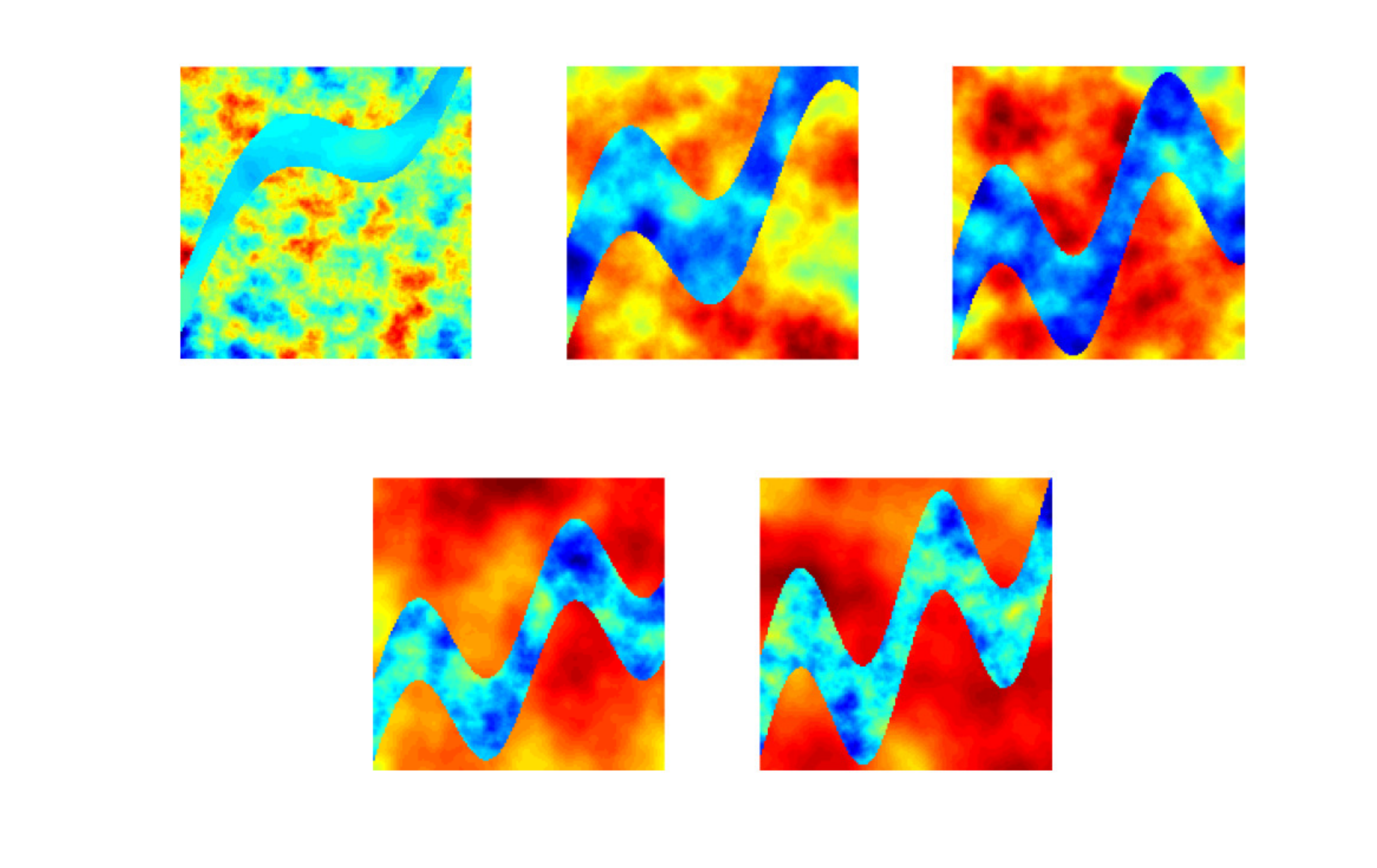}
\caption{Model problem 2. Progression through iterations of non-centered hierarchical method.}
 \label{fig:hier_nc_geo}
\end{figure}

\newpage

\begin{figure}[h!]
\centering
 \includegraphics[scale=0.75]{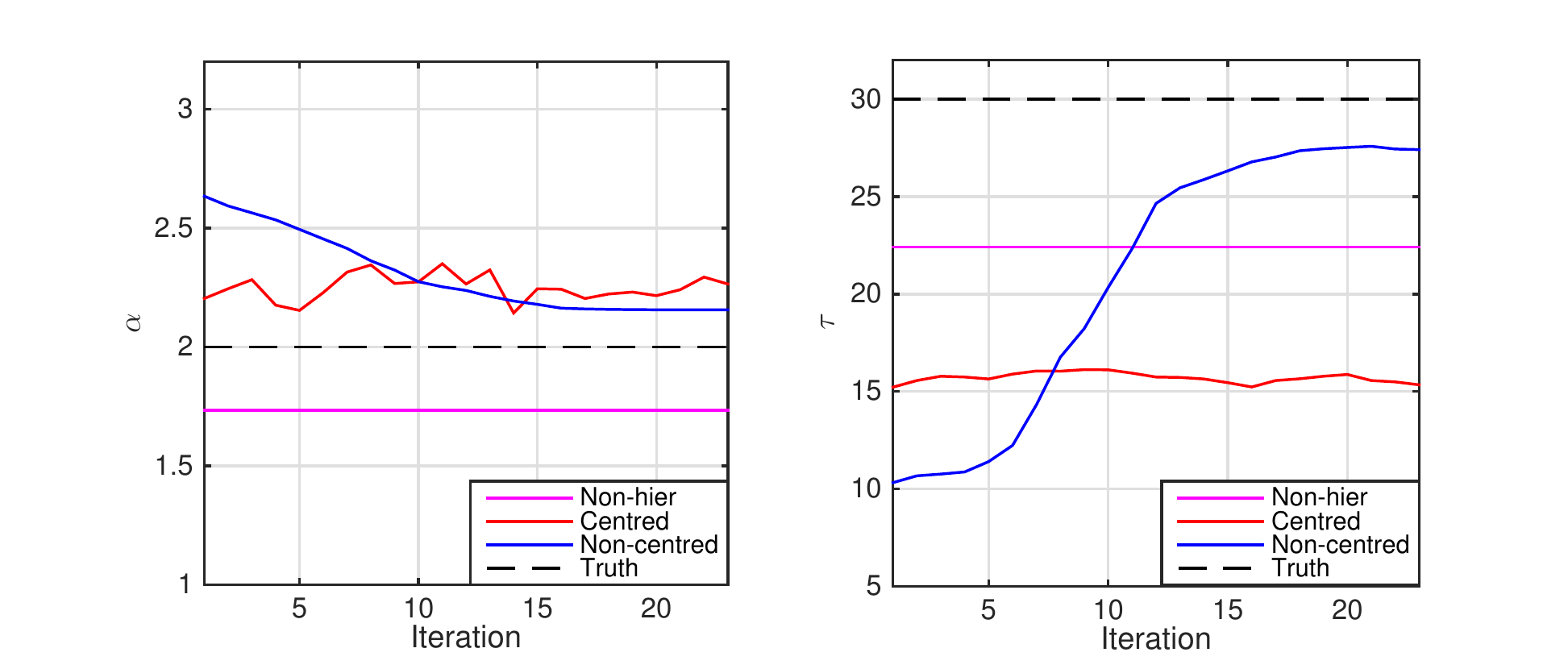}
\caption{Model problem 2. Progression of average value for $\alpha_{1}$ and $\tau_{1}$. }
 \label{fig:new1}
\end{figure}

\newpage

\begin{figure}[h!]
\centering
 \includegraphics[scale=0.75]{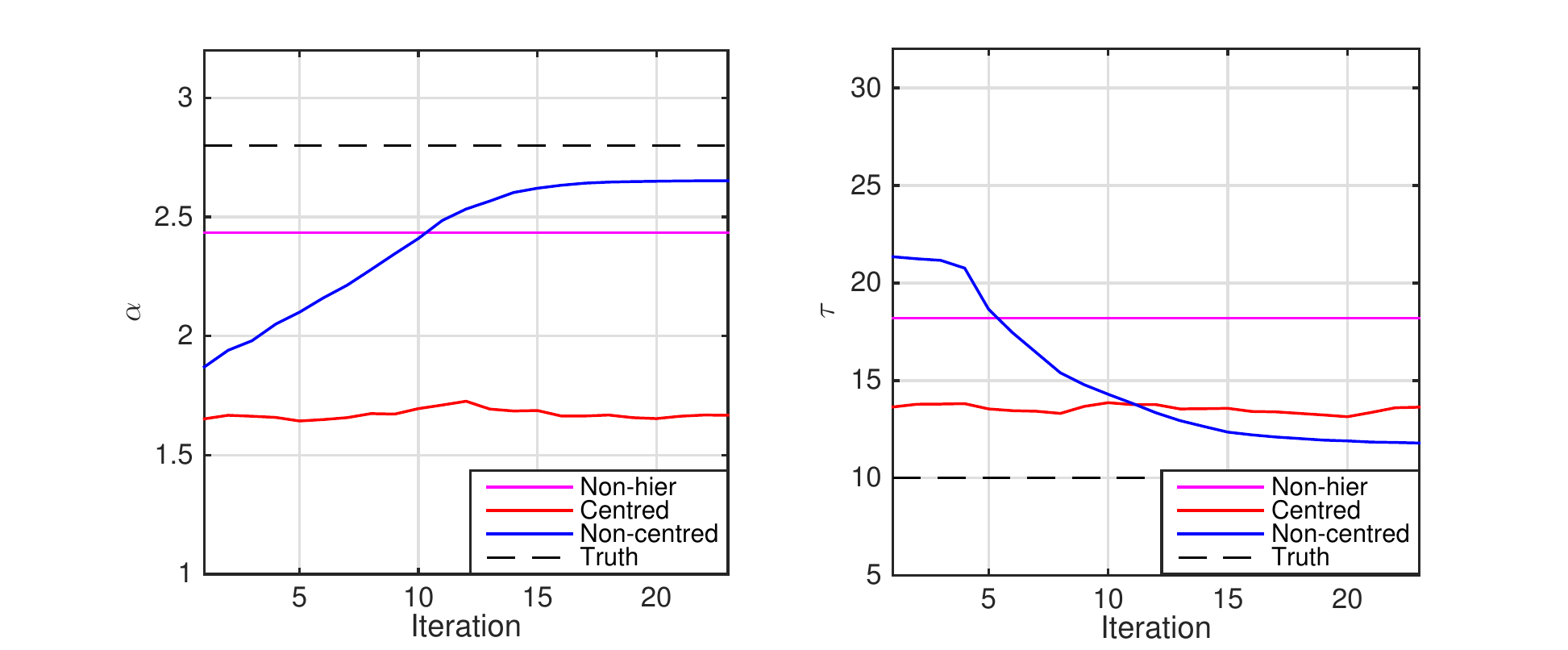}
\caption{Model problem 2. Progression of average value for $\alpha_{2}$ and $\tau_{2}$. }
 \label{fig:new2}
\end{figure}

\begin{figure}[h!]
\centering
 \includegraphics[scale=0.75]{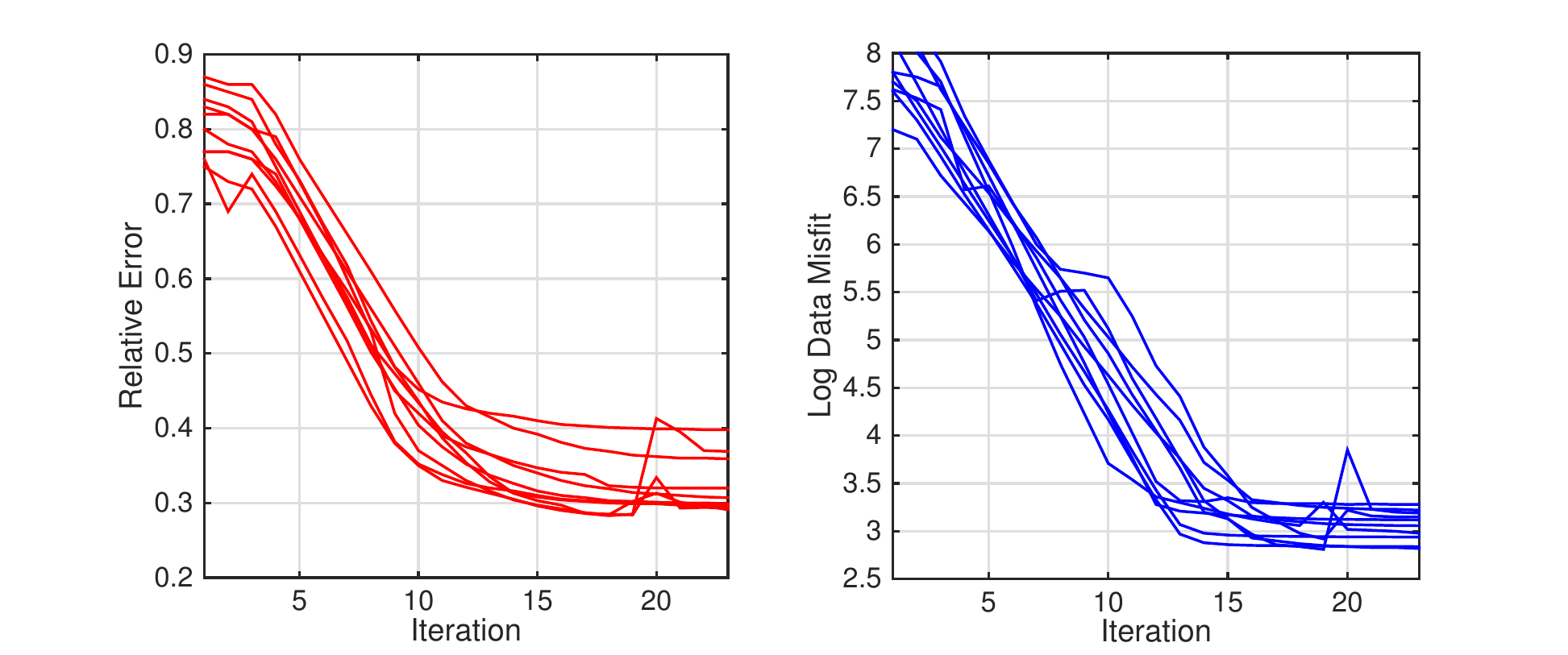}
\caption{Model problem 2. Left: relative error. Right: log-data misfit.}
 \label{fig:channel_error_misfit}
\end{figure}

Figures  \ref{fig:non_hier_geo} -  \ref{fig:hier_nc_geo} are consistent with the previous subsection
in that we notice that hierarchical methods are needed and that non-centring
is necessary to make hierarchical methods perform well. We qualify this
by noting that the geometry is well-learned in all cases, but that the reconstructions of the 
random fields $u_{1}$ and $u_{2}$ inside and outside the geometry are sensitive
to needing non-centered hierarchical representation. Even then the
reconstruction is only accurate in terms of amplitude and length scales
and not pointwise.  

\begin{figure}[h!]
\centering
 \includegraphics[scale=0.67]{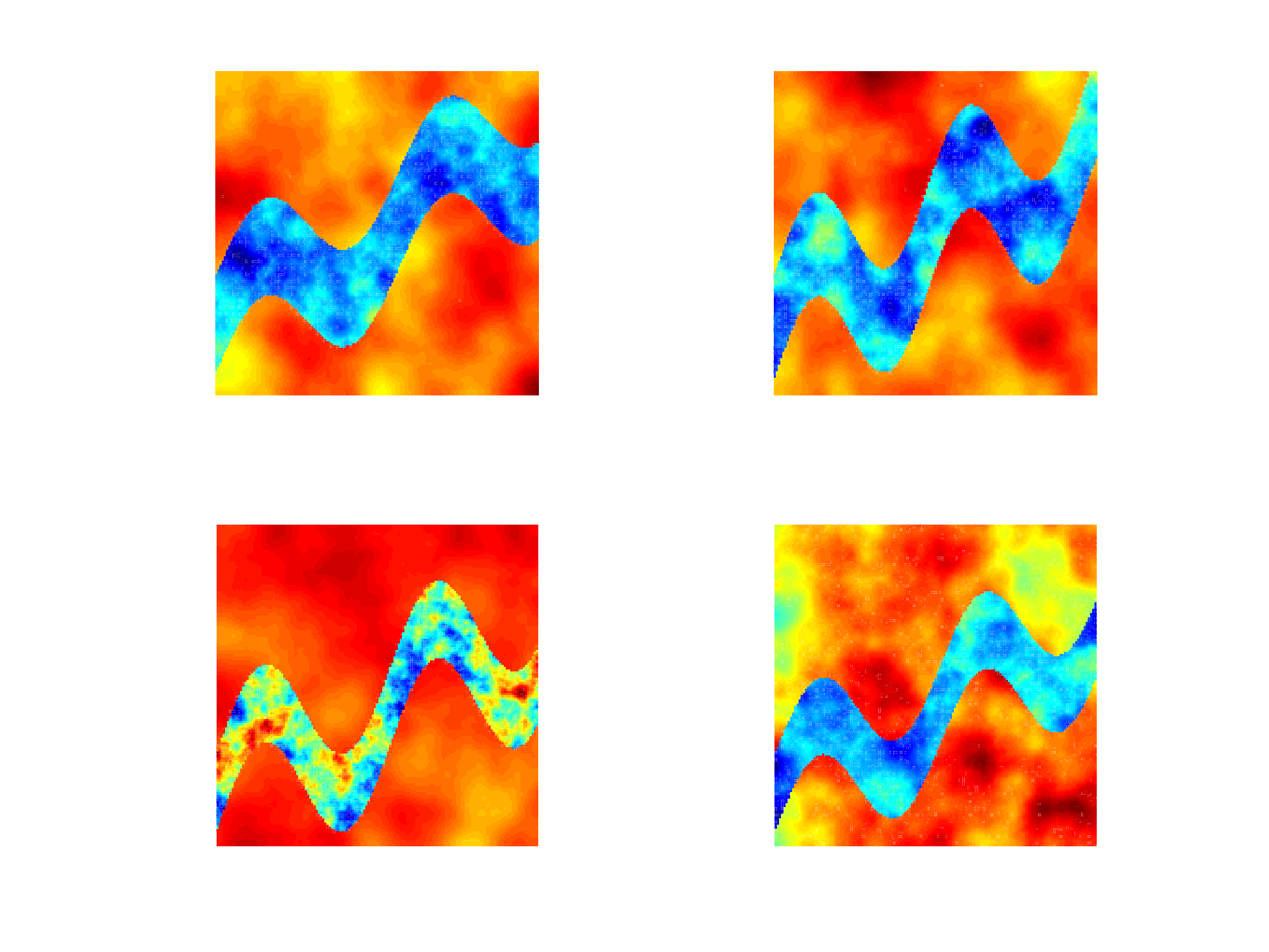}
\caption{Model problem 2. EKI for the final iteration for the non-centered approach from four different initializations. }
 \label{fig:4_final_channel}
\end{figure}

{Figures \ref{fig:new1} and \ref{fig:new2}} give further insight into this, showing how the
smoothness and length scale parameters are learned differently in the
hierarchical, centered and non-centered methods. Again the conclusions are
consistent with the previous subsection.
Figures \ref{fig:channel_error_misfit} and \ref{fig:4_final_channel} 
concentrate on the application of the non-centered approach, using
ten different initializations. The data misfit and 
relative error in the field are shown for all ten cases 
in Figure \ref{fig:channel_error_misfit}; four solution 
estimates are displayed in Figure \ref{fig:4_final_channel}.
The results are similar to those in the previous subsection. However
there is more variability across initializations. This might be
ameliorated by use of a larger ensemble size. 

\newpage

\subsection{Function-valued Hierarchical Parameterization}\label{subsection-4.3}
Our final set of experiments will be based on hierarchical inversion of non-stationary random fields, using Model Problem 3. We consider reconstruction of 
truths that are not drawn from the prior; the prior will be a hierarchical
Gaussian model with spatially varying inverse length scale as {hyperparameter}.
Examples of such truths are ones which contain both rough and smooth features.  We discretize the forward model using a piecewise-linear finite element method (FEM), with a mesh of $h=1/100$. The truth we aim to recover is given by 
\begin{equation}
\label{eq:step}
u^{\dagger}(x)=
\begin{cases}
\exp \bigg(4 - \frac{25}{x(5-x)} \bigg), & x \in (0,5)\\
1, & x \in [7,8]\\
-1, & x \in (8,9]\\
0, & \textrm{otherwise},
\end{cases} \bigskip
\end{equation}
which incorporates both rough and smooth features. Our parameters for the iterative methods are identical to those used in subsection  \ref{subsection-4.1}. 
Because the results of subsection \ref{subsection-4.1} and \ref{subsection-4.2}
clearly demonstrate the need for non-centring in hierarchical methods, 
we do not include results for the centered hierarchical approach here; 
we compare non-hierarchical methods with the use of 
non-centered hierarchical methods with both Cauchy and Gaussian random 
fields as priors on the {hyperparameter} $v$.

\begin{figure}[h!]
\centering
 \includegraphics[scale=0.63]{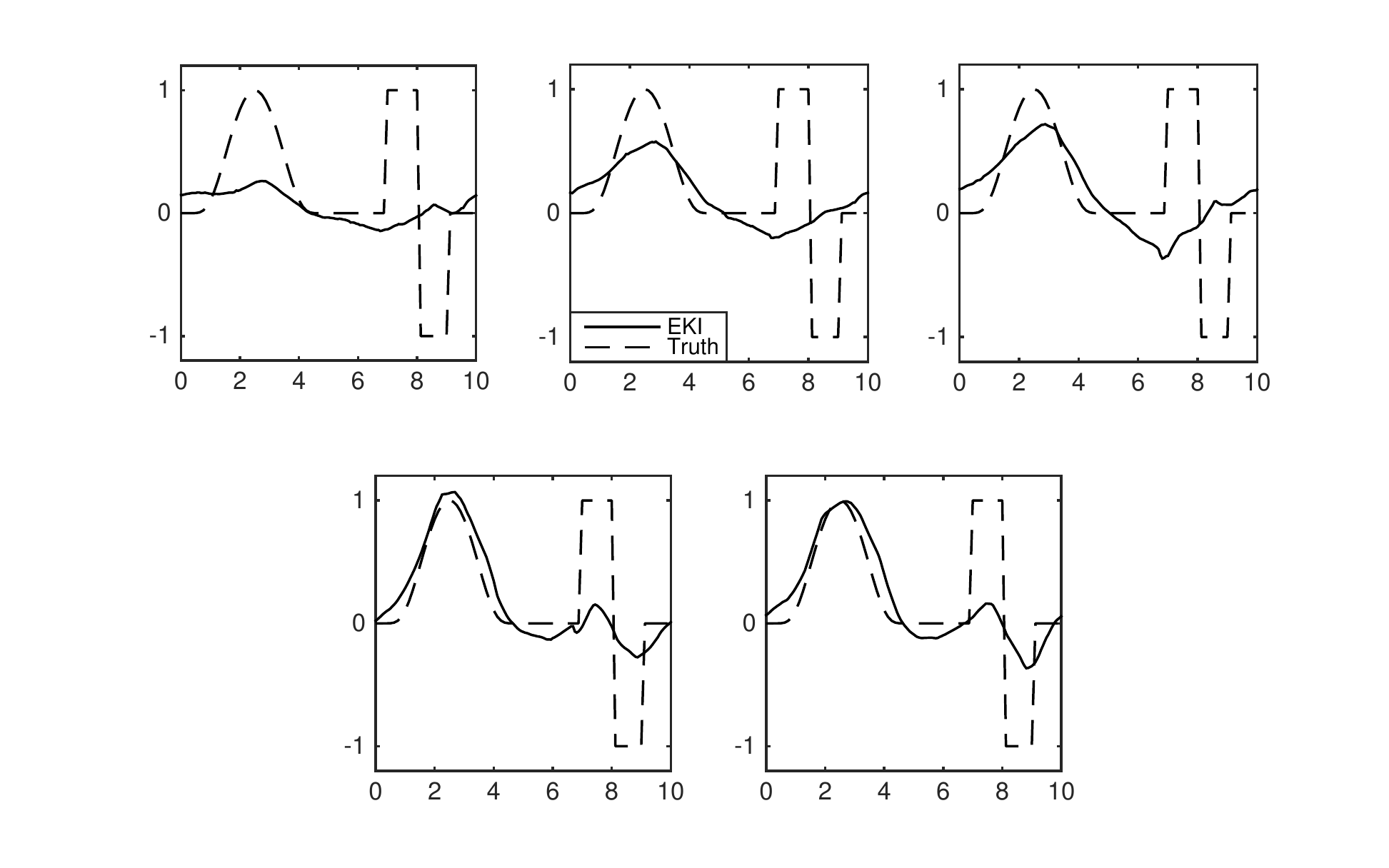}
\caption{Model problem 3. Progression through iterations of non-hierarchical method.}
 \label{fig:none}
\end{figure}

\begin{figure}[h!]
\centering
 \includegraphics[scale=0.63]{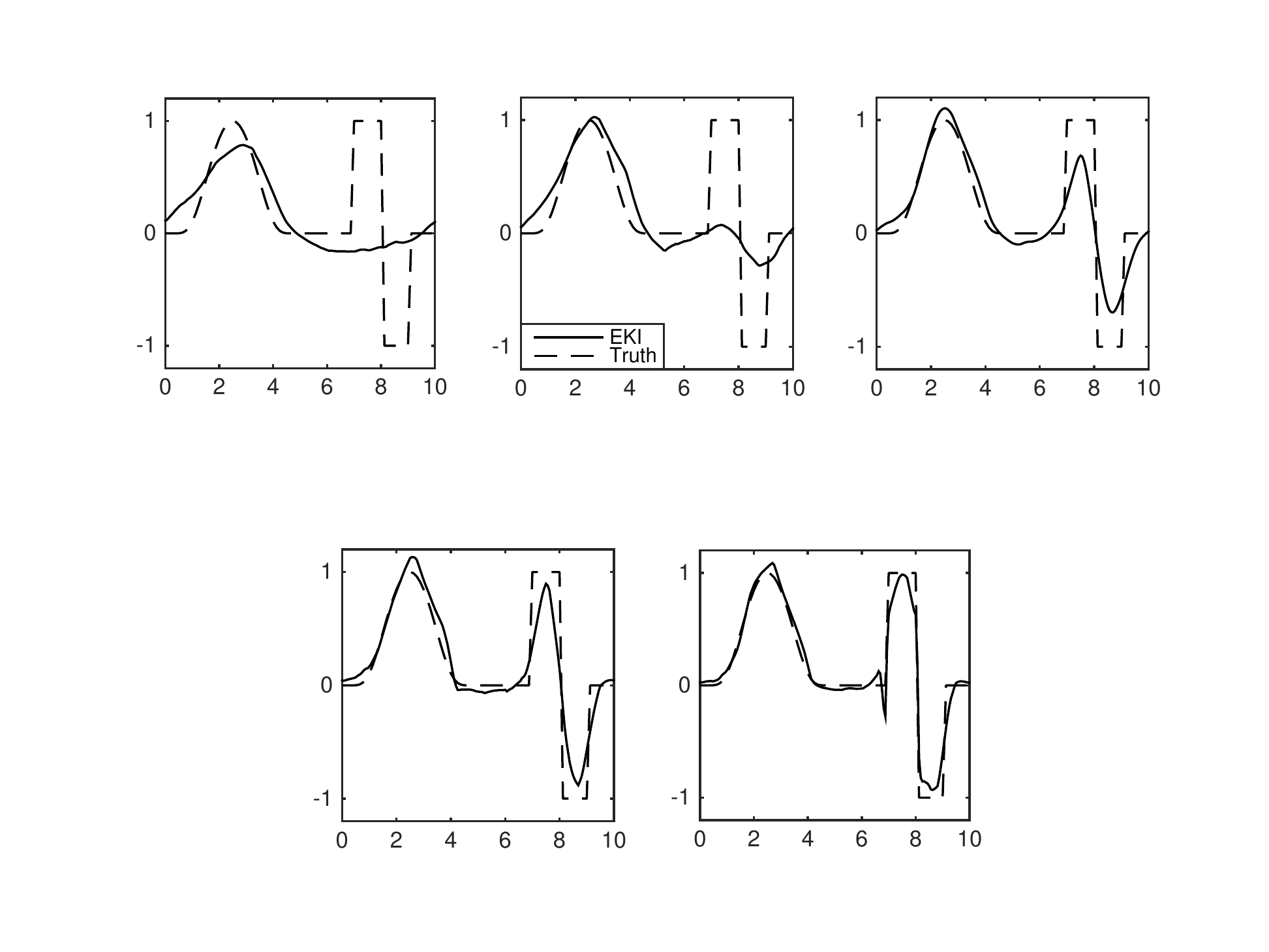}
\caption{Model problem 3. Progression through iterations with hierarchical Gaussian random field.}
 \label{fig:Gauss}
\end{figure}

\newpage

\begin{figure}[h!]
\centering
 \includegraphics[scale=0.63]{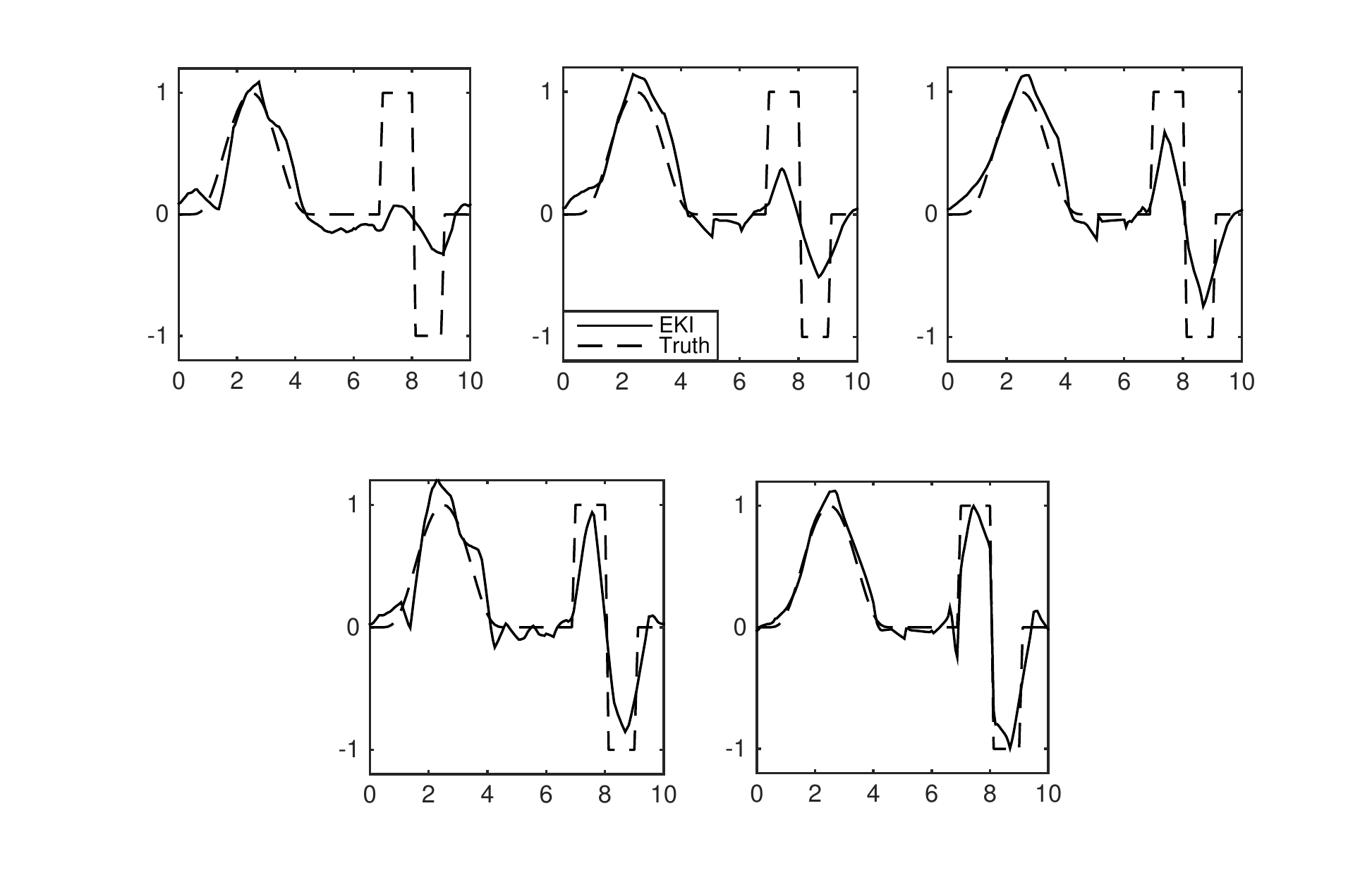}
\caption{Model problem 3. Progression through iterations with hierarchical Cauchy random field.}
 \label{fig:Cauchy}
\end{figure}

\begin{figure}[h!]
\centering
 \includegraphics[scale=0.75]{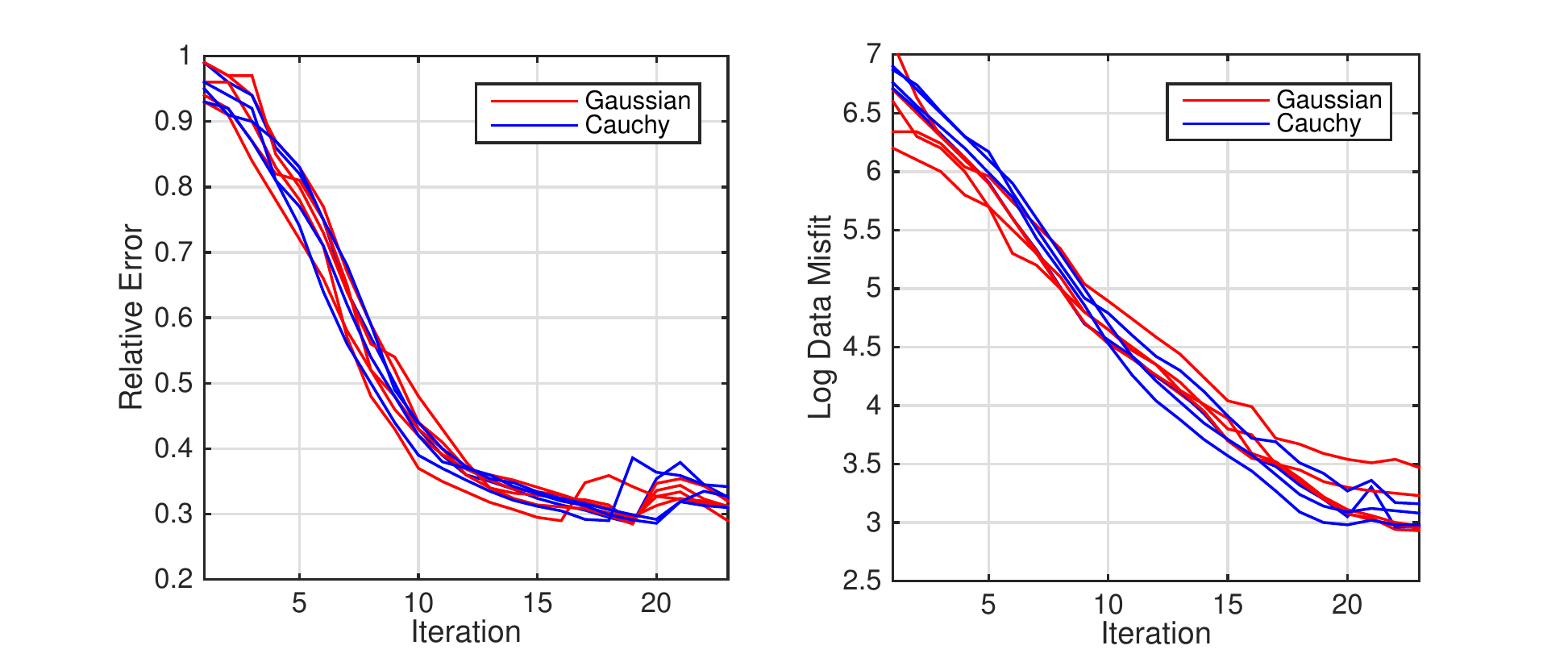}
\caption{Model problem 3. Left: relative error. Right: log-data misfit.}
 \label{fig:field_error_misfit}
\end{figure}

\begin{figure}[h!]
\centering
 \includegraphics[width=\linewidth]{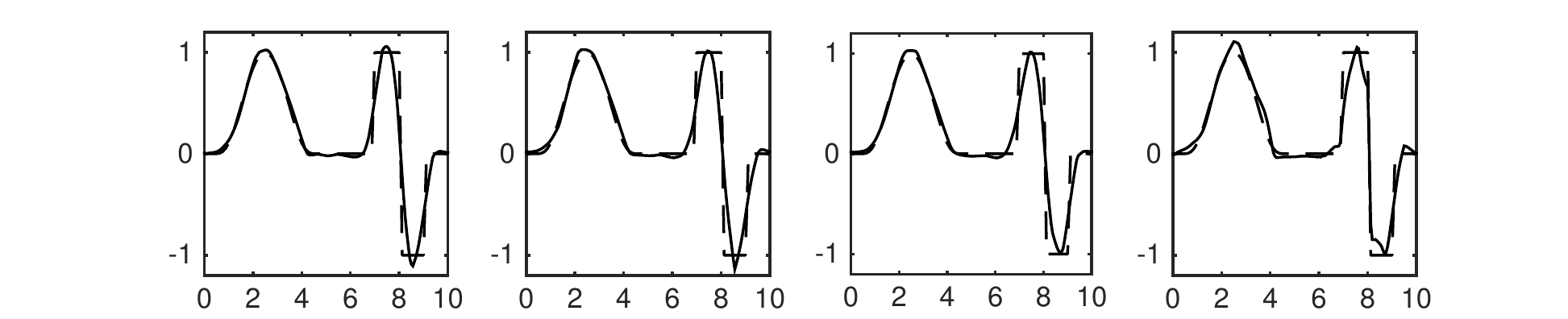}
\caption{Model problem 3.  EKI for the final iteration for the Gaussian hierarchical method from four different initializations.}
 \label{fig:gauss_nc}
\end{figure}

\begin{figure}[h!]
\centering
 \includegraphics[width=\linewidth]{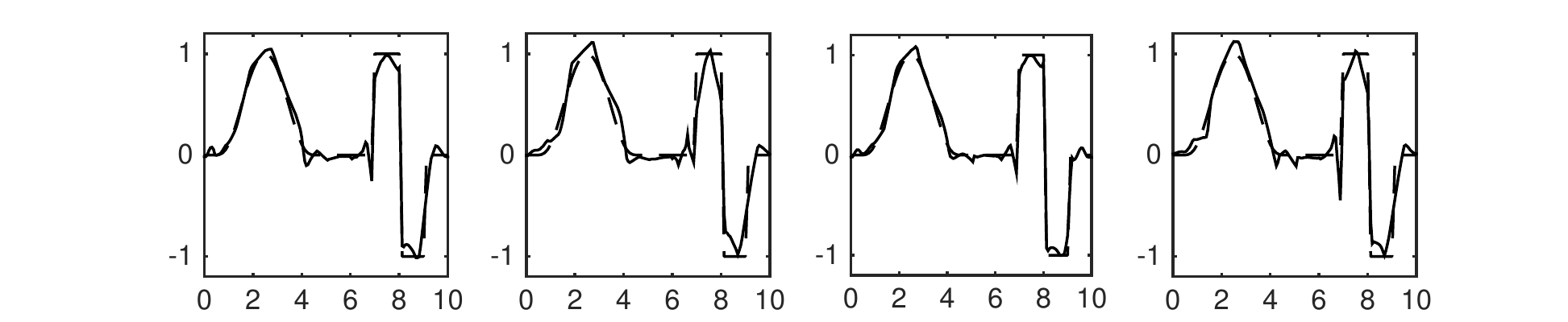}
\caption{Model problem 3.  EKI for the final iteration for the Cauchy hierarchical method from four different initializations.}
 \label{fig:cauchy_nc}
\end{figure}

\newpage

For our Cauchy density function \eqref{eq:cp}, we set $a=4$ and $b=d=0$. Our length scale for the non-hierarchical method will be based on a Gaussian random field \eqref{eq:gp}, similarly with the non-centered Gaussian approach. Our comparison of each approach is provided in Figures \ref{fig:none} - \ref{fig:Cauchy}. We notice that the non-hierarchical method struggles to reconstruct the truth,
and in particular the piecewise constant part of it with
discontinuities. In contrast both non-centered approaches perform
well.  
The effectiveness of the non-centered approaches are highlighted in
Figures \ref{fig:field_error_misfit} - \ref{fig:cauchy_nc}. The first shows the data-misfit and relative error over ten realizations,
and the second shows four reconstructions chosen from these ten at random.
Substantial robustness to the choice of realization is clear.
Note however that for a multi-dimensional problem inversion the
variability with respect to initialization may be more significant,
as in the previous two subsections.

\section{Conclusion \& Discussion} \label{section-5}
In this paper we have considered several forms of parameterizations for 
EKI. In particular our main contribution has been to
highlight the potential for the use of hierarchical techniques from 
computational statistics, and the use of geometric parameterizations, such
as the level set method. Our perspective on EKI is
that it forms a derivative free optimizer and we do not evaluate it from
the perspective of uncertainty quantification. However the hierarchical and
level set ideas are motivated by Bayesian formulations of inverse problems.
We have shown that our parameterizations do indeed lead to better 
reconstructions of the truth on a variety of model problems including
groundwater flow, EIT and source inversion.  
There is very little analysis of EKI, especially
in the fixed, small, ensemble size setting where it is most powerful.
Existing work in this direction may be found in  \cite{Dirk,Claudia}; it
would be of interest to extend these analyses to the parameterizations
introduced here. Furthermore, from a practical perspective, it would
be interesting to extend the deployment of the methods introduced here
to the study of further applications.

\bigskip
\noindent{\bf Acknowledgements.} The authors are grateful to {M. M. Dunlop},
O.\ Papaspiliopoulos and G.\ O.\ Roberts for helpful discussions about centering and non-centering parameterizations. The research  of AMS  was partially supported by the EPSRC programme grant EQUIP and by AFOSR Grant FA9550-17-1-0185. LR was supported by the EPSRC programme grant EQUIP. NKC was partially supported by the EPSRC MASDOC Graduate Training Program and by Premier Oil.

{}


\section*{References}
\begin{thebibliography}{8}
\bibitem{Aan} Aanonsen S I,  Naevdal G,  Oliver D S,  Reynolds A C  and Valles B 2009 The Ensemble Kalman Filter in Reservoir Engineering-a Review \emph{SPE J.} \textbf{14} 393-412.
\bibitem{Anderson}  Anderson J L 2001 An ensemble adjustment Kalman filter for data assimilation. \emph{Monthly Weather Review} \textbf{129} 2884.
\bibitem{uses}  Adler A and Lionheart W R B 2006 Uses and abuses of EIDORS: An extensible software base for EIT \emph{Physiol Meas} \textbf{27} 25-42.
\bibitem{Sergios} Agapiou S,  Bardsley J M, Papaspiliopoulos O and Stuart A M 2014 Analysis of the Gibbs sampler for hierarchical inverse problems. \emph{SIAM Journal on Uncertainty Quantification} \textbf{1} 511-544.
\bibitem{IM}  Bakushinsedky A B and Kokurin M Y 2004 Iterative Methods for Approximate Solution of Inverse Problems \emph{Mathematics and its applications}, Springer.
\bibitem{GN}  Bauer F, Hohage T and Munk A 2009 Iteratively Regularized Gauss-Newton Method for Nonlinear Inverse Problems with Random Noise \emph{SIAM J. Numer. Anal} \textbf{47} 1827-1846.
\bibitem{Dirk} Bl\"{o}mker D, Schillings C and Wacker P 2017 A strongly convergent numerical scheme from EnKF continuum analysis arXiv:1703.06767.
\bibitem{EIT}  Borcea L 2002 Electrical Impedance Tomography \emph{Inverse Problems} \textbf{18} 99-136.
\bibitem{Burger1} Burger M 2002   A  framework  for  the  construction  of  level  set  methods  for  shape  optimization  and reconstruction \emph{Interfaces and Free Boundaries} \textbf{5} 301-329.
\bibitem{Burger2}  Burger M and Osher S 2005  A survey on level set methods for inverse problems and optimal design \emph{Eur. J. Appl. Math} \textbf{16} 263-301.
\bibitem{BC}  Carrera J and Neuman S P 1986 Estimation of aquifer parameters under transient and steady state conditions: application to synthetic and field data \emph{Water Resources Research} \textbf{22} 228-242.
\bibitem{NKC2}Chada N K 2018 {Ensemble Based Methods for Geometric Inverse Problems} \emph{PhD Dissertation}.
\bibitem{BIP} Dashti M Stuart A M 2017 The Bayesian Approach To Inverse Problems \emph{Handbook of Uncertainty Quantification} 311-428.
\bibitem{VC1} {Chen V}, Dunlop M M, Papaspiliopoulos O and Stuart A M 2017 
Robust MCMC Sampling with Non-Gaussian and Hierarchical Priors in High Dimensions,
\emph{In preparation}.
\bibitem{MD} Dunlop M M and Stuart A M 2016 The Bayesian formulation of EIT \emph{Inverse Problems and Imaging} \textbf{10} 1007-1036.
\bibitem{MD2} Dunlop M M, Iglesias M A and Stuart A M 2016 Hierarchical Bayesian Level Set Inversion \emph{Statistics and Computing}.
\bibitem{Hier} Emerick A A 2016 Towards a hierarchical parametrisation to address prior uncertainty in ensemble-based data assimilation \emph{Computational Geosciences} \textbf{20} 35-47.
\bibitem{ernst2015analysis}
Ernst, O and Sprungk, B and Starkloff, H. 2015 Analysis of the ensemble and
  polynomial chaos {K}alman filters in {B}ayesian inverse problems, 
\emph{SIAM-ASA JUQ} \textbf{3} 823--851.
\bibitem{OG} Evensen 1994 G Sequential data assimilation with a nonlinear quasi-geostrophic model using Monte Carlo methods to forecast error statistics \emph{J. Geophysical Research-All Series} \textbf{99} 10-10.
\bibitem{EnKF} Evensen G 2009 \emph{Data Assimilation: The Ensemble Kalman Filter} Springer.
\bibitem{VL} Evensen G and Van Leeuwen P J 1996 Assimilation of geosat altimeter data for the agulhas current using the ensemble Kalman filter with a quasi-geostrophic model \emph{Monthly Weather Review} \textbf{128} 85-96.
\bibitem{Hanke} Hanke M 1997 A regularizing Levenberg-Marquardt scheme, with applications to inverse groundwater filtration problems \emph{Inverse Problems} \textbf{13} 79-95.
\bibitem{Weather} Houtekamer P L and Mitchell H L 2001  A sequential ensemble Kalman filter for atmospheric data assimilation \emph{Monthly Weather Review} \textbf{129} 123-137.
\bibitem{reg} Iglesias M A 2016 A regularising iterative ensemble Kalman method for PDE-constrained inverse problems \emph{Inverse Problems} \textbf{32} .
\bibitem{ILS} Iglesias M A, Law K J H and Stuart A M 2013 Evaluation of Gaussian approximations for data assimilation in reservoir models \emph{Computational Geosciences}, Volume \textbf{17} 851-885.
\bibitem{ILS2} Iglesias M A, Law K J H and Stuart A M 2014 Ensemble Kalman methods for inverse problems \emph{Inverse Problems} \textbf{29}.
\bibitem{geo} Iglesias M A, Lin K, and Stuart A M 2014 Well-posed Bayesian geometric inverse problems arising in subsurface flow \textbf{30}.
\bibitem{Lu} Iglesias M A, Lu Y and Stuart A M 2015 A Bayesian Level Set Method for Geometric Inverse Problems \emph{Interfaces and Free Boundary Problems}.
\bibitem{Kaipo} Kaipo J and Somersalo E 2004 \emph{Statistical and Computational Inverse Problems} {Applied Mathematical Sciences}, Springer.
\bibitem{Kalman} Kalman R E 1960 A new approach to linear filtering and prediction problems \emph{Trans ASME (J. Basic Engineering)} {\bf 82}  35-45.
\bibitem{Kantas} Kantas N, Beskos A and Jasra A 2014 Sequential Monte Carlo Methods for High-Dimensional Inverse Problems: A case study for the Navier-Stokes equations \emph{SIAM/ASA J. Uncertain. Quantif} {\bf 2} 464-489.
\bibitem{Li} Li G and Reynolds A C 2009 Iterative ensemble Kalman filters for data assimilation \emph{SPE J} {\bf 14} 496-505.
\bibitem{Majda}  Majda A and Wang W 2006 \emph{Non-linear Dynamics and Statistical Theories for Basic Geophysical Flows} Cambridge University Press.
\bibitem{Markku} Markkanen M, Roininen L, Huttunen J M J and Lasanen S 2015 Cauchy difference priors for edge-preserving Bayesian inversion with an application to X-ray tomography arXiv:1603.06135.
\bibitem{LS} Law K J H and Stuart A M 2012 Evaluating Data Assimilation Algorithms \emph{Mon. Weather Rev} {\bf 140} 37-57.
\bibitem{Liang}  Liang P, Jordan M I and Klein D 2010 Learning Programs: A Hierarchical Bayesian Approach \emph{Conference: Proceedings of the 27th International Conference on Machine Learning}.
\bibitem{Lind} Lindgren F, Rue H and Lindstr\"om J 2011 An explicit link between Gaussian fields and Gaussian Markov random fields: the stochastic partial differential equation approach \textbf{73} 423-498.
\bibitem{Liu} Liu W, Li J and Marzouk Y M 2017 An approximate empirical Bayesian method for large-scale linear-Gaussian inverse problems  arXiv:1705.07646.
\bibitem{SPE}  Myserth I and Omre H 2010 Hierarchical Ensemble Kalman Filter \emph{SPE J} \textbf{15}.
\bibitem{Oliver} Oliver D, Reynolds A C and Liu N 2008 \emph{Inverse Theory for Petroleum Reservoir Characterization and History Matching} Cambridge University Press.
\bibitem{centre1} Papaspiliopoulos O, Roberts G O and Sk\"{o}ld M 2003 Non-centered parameterisations for hierarchical models and data augmentation In \emph{Bayesian Statistics}  7:  Proceedings of the Seventh Valencia International Meeting. \emph{Oxford University Press}.
\bibitem{centre2} Papaspiliopoulos O, Roberts G O and Sk\"{o}ld M 2007 A general framework for the parametrization of hierarchical models \emph{Statistical Science} 59-73.
\bibitem{SMC} Pierre Del Moral A J  and Doucet A 2006 Sequential Monte Carlo samplers \emph{Journal of the Royal Statistical Society} Series B 411-436.
\bibitem{lassi_hier}  Roininen L, Girolami M, Lasanen S and Markkanen M 2017 Hyperpriors for Matérn fields with applications in Bayesian inversion arXiv:1612.02989.
\bibitem{lassi} Roininen L, Huttunen J M J and Lasanen S 2014 Whittle-Mat\'{e}rn priors for Bayesian statistical inversion with applications in electrical impedance tomography \emph{Inverse problems and imaging} \textbf{8} .561-586.
\bibitem{raul} Ruggeri F, Sawlan Z, Scavino M and Tempone R 2016 A Hierarchical Bayesian Setting for an Inverse Problem in Linear Parabolic PDEs with Noisy Boundary Conditions \emph{Bayesian Anal}.
\bibitem{Santosa1996} Santosa, F 1996 { A level-set approach for inverse problems involving obstacles} {\em ESAIM}, 1(January):17--33, 1996.
\bibitem{Claudia} Schillings C and Stuart A M 2017 Analysis of the ensemble Kalman filter for inverse problems.  SIAM J Numerical Analysis {\bf 55}, 1264-1290. 
\bibitem{EIT2} Somersalo E, Cheney M and Isaacson D 1992 Existence and Uniqueness for Electrode Models for Electric Current Computed Tomography \emph{SIAM J. Appl. Math} \textbf{52} 023-1040
\bibitem{Teh} Teh Y W, Jordan M I, Beal M J, and Blei D M 2016 Hierarchical Dirichlet processes \emph{Journal of the American Statistical Association} \textbf{101} 1566-1581.
\bibitem{SPE-2} Tsyrulnikov M and Rakitko A 2016 Hierarchical Bayes Ensemble Kalman Filter  \emph{Physica D: Nonlinear Phenomena}.
\bibitem{Tar} Tarantola A 1987 \emph{Inverse Problem Theory and Methods for Model Parameter Estimation} {Elsevier}.
\bibitem{ESRKF} Tippett M K 2002 Ensemble Square Root Filters \emph{Monthly Weather Review} \textbf{131}.
\end{thebibliography}
\end{document}